\documentclass{article}
\usepackage{amsmath,amsthm,amsfonts,eucal}

\numberwithin{equation}{section}

\def\ca{{\mathcal A}}
\def\cb{{\mathcal B}}
\def\cc{{\mathcal C}}
\def\cd{{\mathcal D}}

\def\cf{{\mathcal F}}
\def\cg{{\mathcal G}}
\def\ch{{\mathcal H}}

\def\cj{{\mathcal J}}
\def\ck{{\mathcal K}}

\def\cam{{\mathcal M}}
\def\cn{{\mathcal N}}

\def\car{{\mathcal R}}

\def\cu{{\mathcal U}}

\def\bg{{\mathbb G}}

\def\bn{{\mathbb N}}
\def\br{{\mathbb R}}
\def\bz{{\mathbb Z}}

\def\a{\alpha}
\def\b{\beta}
\def\g{\gamma}        \def\G{\Gamma}
\def\d{\delta}        \def\D{\Delta}
\def\eps{\varepsilon}

\def\th{\vartheta}

\def\l{\lambda}       \def\La{\Lambda}
\def\m{\mu}
\def\n{\nu}

\def\r{\rho}

\def\t{\tau}

\def\f{\varphi}
\def\c{\chi}
\def\om{\omega}        \def\O{\Omega}

\def\imply{\Rightarrow}
\def\coimply{\Leftarrow}

\def\ov{\overline}
\def\e#1{{\rm e}^{#1}}
\def\itm#1{\item[$(#1)$]}

\def\limr{\lim_{r\to\infty}}
\def\limR{\lim_{R\to\infty}}
\def\lsup{\limsup_{R\to\infty}}
\def\linf{\liminf_{R\to\infty}}

\DeclareMathOperator{\inj}{inj}
\DeclareMathOperator{\supp}{supp}

\def\ad#1{d_\infty(#1)}
\def\md#1{d_0(#1)}

\def\dE{{\partial E}}

\def\dec{\searrow}
\def\aff{\hat\in}

\DeclareMathOperator{\Lim}{Lim}

\def\lo{\Lim_{\om}}
\def\ar{\ca^{\car}}

\def\bg{C$^{\infty}$-bounded geometry }
\def\mt{\ov \cam}
\def\sa{{\text{sa}}}
\def\tordim{{\mathfrak{tordim}}}

\newtheorem{Thm}{Theorem}[section]
\newtheorem{Cor}[Thm]{Corollary}
\newtheorem{Prop}[Thm]{Proposition}
\newtheorem{Lemma}[Thm]{Lemma}
\theoremstyle{definition}
\newtheorem{Dfn}[Thm]{Definition}
\newtheorem{exmp}[Thm]{Example}

\theoremstyle{remark}
\newtheorem{rem}[Thm]{Remark} 
\newtheorem{ack}{Acknowledgement} 
 \title{\huge Novikov-Shubin invariants and asymptotic 
 dimensions for open manifolds}
 \author{Daniele Guido, Tommaso Isola\\
 Dipartimento di Matematica,\\ Universit\`a di Roma ``Tor Vergata'',\\ 
 I--00133 Roma, Italy.\\
{\tt e-mail: guido@mat.uniroma2.it, isola@mat.uniroma2.it}}
\date{July 28, 1998}
\begin{document}
\maketitle
\markboth{Asymptotic dimension package}
{Asymptotic dimension and Novikov-Shubin invariants}
\renewcommand{\sectionmark}[1]{}
\bigskip

\begin{abstract} 
	The Novikov-Shubin numbers are defined for open manifolds with 
	bounded geometry, the $\Gamma$-trace of Atiyah being replaced by a 
	semicontinuous semifinite trace on the C$^*$-algebra of almost 
	local operators.  It is proved that they are invariant under 
	quasi-isometries and, making use of the theory of singular traces 
	for C$^*$-algebras developed in \cite{GI4}, they are interpreted 
	as asymptotic dimensions since, in analogy with what happens in 
	Connes' noncommutative geometry, they indicate which power of the 
	Laplacian gives rise to a singular trace.  Therefore, as in 
	geometric measure theory, these numbers furnish the order of 
	infinitesimal giving rise to a non trivial measure.  The 
	dimensional interpretation is strenghtened in the case of the 0-th 
	Novikov-Shubin invariant, which is shown to coincide, under 
	suitable geometric conditions, with the asymptotic counterpart of 
	the box dimension of a metric space.  Since this asymptotic 
	dimension coincides with the polynomial growth of a discrete 
	group, the previous equality generalises a result by Varopoulos 
	\cite{Varopoulos1} for covering manifolds.
\end{abstract}

\newpage
\setcounter{section}{-1}
\section{Introduction.}\label{sec:intro}

 In a celebrated paper \cite{Atiyah}, Atiyah  observed that on covering 
 manifolds $\G\to M \to X$, a trace on $\Gamma$-periodic operators may 
 be defined, called $\Gamma$-trace, with respect to which the Laplace 
 operator has compact resolvent.  Replacing the usual trace with the 
 $\Gamma$-trace, he defined the $L^2$-Betti numbers and proved an 
 index theorem for covering manifolds.

 Motivated by this paper, Novikov and Shubin \cite{NS1} observed that, 
 since for noncompact manifolds the spectrum of the Laplacian is not 
 discrete, new global spectral invariants can be defined, which measure the 
 density near zero of the spectrum.
 
 Novikov-Shubin invariants and other L$^{2}$-invariants have been a 
 very active research field since then, and the interested reader is 
 referred to \cite{BFK,CFM,Farber,Lott1,Lueck,Lueck1} for recent developments 
 and extensive bibliographies.

 Also based on Atiyah's paper, Roe defined $L^2$-Betti numbers for 
 open manifolds \cite{Roe1} by replacing the $\Gamma$-trace of 
 Atiyah with a trace on a subalgebra of $L^2(M)$, and showed their 
 invariance under quasi-isometries \cite{RoeBetti}.
 
 In this paper, inspired by Roe \cite{RoeCoarse,Roe1}, we define the 
 C$^*$-algebra of almost local operators and a semicontinuous 
 semifinite trace on it, and use this trace to define Novikov-Shubin 
 numbers for open manifolds, proving that they are invariant under 
 quasi-isometries.

 The second part of this paper is concerned with a dimensional 
 interpretation of these numbers.

 As it is known, a general understanding of the geometric meaning of 
 the Novikov-Shubin invariants is still lacking.  To this end, the 
 definition of these numbers in the case of open manifolds corresponds 
 to the idea that interpreting them as global invariants of an open 
 manifold, rather then as homotopy invariants of a compact one, some 
 aspects can be better understood.

 The asymptotic character of these numbers is manifest in two parts of 
 their construction.

 On the one hand, the trace used to define these numbers is a large 
 scale trace, since, as observed by Roe \cite{Roe1}, it is given by an 
 average on the group, in the case of coverings, and by an average on 
 the exhaustion, in the case of open manifolds.

 On the other hand these numbers are defined in terms of the low 
 frequency behaviour of the $p$-Laplacians, or the large time 
 behaviour of the $p$-heat kernel.

 In this respect, they are the large scale counterpart of the spectral 
 dimension, namely of the dimension as it is recovered by the Weyl 
 asymptotics.  Indeed the inverse of the dimension of a manifold 
 coincides with the order of infinitesimal of $\Delta^{-1/2}$, namely 
 with the order of infinitesimal of the eigenvalue sequence 
 $\mu_n(\Delta^{-1/2})$ when $n\to\infty$.  The $p$-th Novikov Shubin 
 invariant $\alpha_p$, instead, coincides with the inverse of the 
 order of infinite, when $t\to0$, of the generalized eigenvalue 
 sequence $\mu_t(\Delta_p^{-1/2})$.

 A fundamental observation of Connes is that integration on a compact
 Riemannian manifold may be reconstructed by making use of the
 logarithmic trace and the Weyl asymptotics.  Indeed if the resolvent
 of the Dirac operator is compact and is an infinitesimal of order
 $1/d$, then $|D|^{-d}$ is (logarithmically) traceable, and the
 corresponding singular trace reconstructs the integration on the
 manifold.

 This remark goes in the direction of a noncommutative geometric 
 measure theory.  In fact, while in geometric measure theory the 
 dimension is the unique exponent to give to the radius of a ball in 
 order to obtain, using Hausdorff procedure, a (possibly) non trivial 
 measure, in noncommutative geometry the dimension may be defined as 
 the exponent to give to the resolvent of the Dirac operator in order 
 to obtain a non-trivial singular trace, hence a non-trivial measure 
 on the given space.

 Also in this respect the Novikov-Shubin numbers may well be 
 considered asymptotic spectral dimensions.  In \cite{GI1} a new type 
 of singular traces, for continuous semifinite von~Neumann algebras, 
 were introduced, the so called singular traces at $0$, which measure 
 the divergence at $0$ of the generalized eigenvalue function $\mu(t)$ 
 introduced by Fack and Kosaki \cite{FK}.  Such traces were then 
 defined also in the case of C$^*$-algebras via the noncommutative 
 Riemann integration \cite{GI4}.  We show, in analogy with the local 
 results, that the operator $\Delta_p^{-1/2}$, raised to the power 
 $\alpha_p$, is singularly traceable, hence a singular trace is 
 naturally attached to these asymptotic spectral dimensions.

\medskip

 The only Novikov-Shubin number for which a clear geometric
 interpretation has been given is $\alpha_0$, in fact Lott noted in
 \cite{Lott} that a result of Varopoulos \cite{Varopoulos1}
 immediately implies the equality of $\alpha_0$ with the growth of the
 fundamental group in the case of covering manifolds.  We prove a
 generalization of this result in the case of open manifolds with
 bounded geometry and satisfying an isoperimetric inequality
 introduced by Grigor'yan \cite{Grigoryan94}. As a consequence, we 
 compute the range of $\a_{0}$ for such manifolds to be $[1,\infty)$.

 First, we associate a number to any metric space, which we
 call asymptotic dimension since it is the global analogue of the box
 dimension defined by Kolmogorov and Tihomirov \cite{KT}, and show
 that it is invariant under rough isometries.  This shows in
 particular that the asymptotic dimension of a covering manifold
 coincides with the growth of the fundamental group.  
 
 Then, for the mentioned class of open manifolds, we show that such 
 asymptotic dimension coincides with the $0$-th Novikov-Shubin number.  
 This result, besides strengthening our dimensional interpretation of 
 $\alpha_0$, shows a stronger invariance property for it.  Indeed the 
 Novikov-Shubin numbers depend on a chosen regular exhaustion of the 
 manifold, hence on its geometry in the large.  When the mentioned 
 isoperimetric inequality holds, the volume growth of the manifold is 
 subexponential, and in this case Roe proved \cite{Roe1} that there is 
 a regular exhaustion given by balls with a fixed center.  With this 
 natural choice, $\alpha_0$ is invariant under rough isometries.

 From the technical point of view, a large part of this paper deals
 with the problem of defining a semicontinuous semifinite trace on the
 C$^*$-algebra of almost local operators, i.e. on the norm closure of
 the operators with finite propagation. Such a trace depends on the
 geometry in the large of the manifold or, more precisely, on the
 exhaustion $\ck$.  The semicontinuity and semifiniteness properties
 allow us to use the theory of noncommutative Riemann integration
 developed in \cite{GI4}, hence to extend the trace to a bigger
 algebra, containing many projections, and eventually to a bimodule of
 unbounded operators affiliated to it.  Then we can associate a
 positive non-increasing function, the generalized eigenvalue function
 $\m_A$, with any operator $A$ in the bimodule, and define the
 asymptotic spectral dimension of $(M,\ck,\D_p)$ as the inverse of the
 ``order of infinite'' of $\mu_{\Delta_p^{-1/2}}(t)$, when $t\to0$,
 showing that, on the one hand, it coincides with the Novikov-Shubin
 number $\alpha_p$, and, on the other hand, it produces a noncommutative
 integration procedure, i.e. a singular trace.

 The semicontinuity and semifiniteness properties of our trace 
 constitute a technical semplification in the study of Novikov-Shubin 
 numbers, and are crucial for the possibility of defining type 
 II$_{1}$ singular traces for C$^{*}$-algebras.  Our proof of the 
 quasi-isometry invariance of the Novikov-Shubin numbers parallels the 
 corresponding invariance for Betti numbers proved in \cite{RoeBetti}.

 The definition of singular traces at $0$ for C$^*$-algebras is
 contained in \cite{GI4}, and is briefly described in section
 \ref{subsec:singtrac}; in the same subsection the singular
 traceability of $\Delta_p^{-\alpha_p}$ is proved, making use of a
 general result contained in \cite{GI5}.  Let us remark that we do not
 need $\mu_{\Delta_p^{-1/2}}(t)$ to have an exact polinomial
 behaviour, when $t\to0$, but only to have a polinomial bound from
 below.  The inverse of $\alpha_p$ then coincides with the exponent of
 the ``optimal'' bound from below or, more precisely, with the
 supremum of the $\om$ such that $t^{-\om}$ is a lower bound for
 sufficiently small $t$.  The existence of a polynomial bound from
 below, namely the positivity of $\alpha_p$, guarantees the singular
 traceability of $\Delta_p^{-\alpha_p/2}$ regardless of the fact that
 such operator belongs to $L^1(\tau)$ or not, in particular the
 singular trace is not necessarily the logarithmic trace introduced 
 by Dixmier \cite{Dixmier}.  This fact corresponds
 to the idea that the logarithmic behaviour has to be expected only in
 the regular cases, such as smooth manifolds, which are locally 
 regular, or covering manifolds, which are regular at large scale, but may
 fail in the general case.

 The definition of asymptotic dimension for metric spaces contained in
 subsection \ref{subsec:asympdim} is obtained from the definition of
 dimension given by Kolmogorov and Tihomirov (often called box
 dimension), simply replacing limits to $0$ with limits to $\infty$
 and viceversa.  More precisely, if $n(r,R)$ denotes the minimum
 number of balls of radius $r$ necessary to cover a ball of radius $R$
 (and given center), the box dimension is the ``order of infinite'' of
 $n(r,R)$ when $r\to0$ (with $R$ fixed, and often independently of
 $R$), whereas the asymptotic dimension is the ``order of infinite''
 of $n(r,R)$ when $R\to\infty$ (with $r$ fixed, and often
 independently of $r$).

 Our asymptotic dimension enjoys all the formal properties of a
 dimension, plus the invariance under rough isometries, which
 corresponds to its large scale character.

 The equality between the asymptotic dimension and the $0$-th
 Novikov-Shubin number is based on the strict relation between the
 asymptotics of the heat kernel $H_{0}(t,x,x)$ and the asymptotics of
 the volume of a ball of radius $\sqrt{t}$, for $t\to\infty$.  Such
 relation is known to hold for open manifolds with some kind of
 polinomial growth, such as manifolds with positive Ricci curvature
 (cf.  e.g. \cite{Daviesbook}) or manifolds satisfying the
 isoperimetric inequality of Grigor'yan \cite{Grigoryan94}.

 We finally remark that in the definition of box dimension there are
 in principle two possible choices, corresponding to the $\limsup$ or
 to the $\liminf$ procedure in the definition of the order of
 infinite.  But only the $\limsup$ gives rise to the correct behaviour
 for cartesian products, namely the dimension of the product is not
 greater than the sum of the dimensions (cf.  \cite{Pontriagin} and
 \cite{KT}).

 Also in the definition of the Novikov-Shubin numbers two possible 
 choices are available and again the $\limsup$ was chosen since 
 it guarantees the singular traceability property.  It is remarkable 
 that these two independent choices agree, giving rise to the equality 
 between the asymptotic dimension and the $0$-th Novikov-Shubin 
 invariant.
 
 Some of the results contained in the present paper have been 
 announced in several international conferences.  In particular we 
 would like to thank the Erwin Schr\" odinger Institute in Vienna, 
 where this paper was completed, and the organisers of the ``Spectral 
 Geometry Program'' for their kind invitation.

 \section{A trace for open manifolds} \label{sec:trace}
 
 This section is devoted to the construction of a trace on (a suitable
 subalgebra of) the bounded operators on $L^{2}(\La^{p}T^{*}M)$, where
 $M$ is an open manifold of bounded geometry.  The basic idea for this
 construction is due to Roe \cite{Roe1}, and is based on a regular
 exhaustion for the manifold.  We shall regularize this trace, in
 order to get a semicontinuous semifinite trace on the C$^*$-algebra
 of almost local operators.  As observed by Roe, this trace is
 strictly related to the trace constructed by Atiyah \cite{Atiyah} in
 the case of covering manifolds.  It may therefore be used to define
 the Novikov-Shubin invariants for open manifolds, as we do in
 subsection \ref{subsec:laplacian}.

 \subsection{Open manifolds of bounded geometry} 
 \label{subsec:openmanifold}
 
 In this subsection we give some preliminary results on open manifolds 
 of bounded geometry that are needed in the sequel.

 Several definitions of bounded geometry for an open manifold (i.e. a 
 noncompact complete Riemannian manifold) are usually considered.  
 They all require some uniform bound (either from above or from below) 
 on some geometric objects, such as: injectivity radius, sectional 
 curvature, Ricci curvature, Riemann curvature tensor etc.  (For all 
 unexplained notions see e.g. Chavel's book \cite{Chavel}).

 In this paper the following form is used, but see 
 \cite{ChavelFeldman} and references therein for a different approach.  

 \begin{Dfn}\label{2.1.1} 
	 Let $(M,g)$ be a complete Riemannian manifold.  
	 We say that $M$ has C$^{\infty}$-bounded geometry if it has 
	 positive injectivity radius, and the curvature tensor is bounded, 
	 together with all its covariant derivatives.
 \end{Dfn}
 
 \begin{Lemma}\label{2.1.2} 
	 Let $M$ be an $n$-dimensional complete Riemannian manifold with 
	 positive injectivity radius, sectional curvature bounded from 
	 above by some constant $c_1$, and Ricci curvature bounded from 
	 below by $(n-1)c_2 g$, in particular M could have 
	 C$^\infty$-bounded geometry.  Then there are real functions 
	 $\b_1,\ \b_2$ s.t. \itm{i} for all $x\in M$, $r>0$,
	 $$
	 0<\b_1(r)\leq vol(B(x,r)) \leq \b_2(r),
	 $$
	 \itm{ii} $\lim_{r\to0}\frac{\b_2(r)}{\b_1(r)} =1$. 
 \end{Lemma}
 \begin{proof}
	 $(i)$ We can assume $c_2<0<c_1$ without loss of generality.  
	 Then, denoting with $V_\d(r)$ the volume of a ball of radius $r$ 
	 in a manifold of constant sectional curvature equal to $\d$, we 
	 can set $\b_1(r) := V_{c_1}(r\wedge r_0)$, and $\b_2:= 
	 V_{c_2}(r)$, where $r_0:= \min\{ \inj(M), \frac{\pi}{\sqrt{c_1}} 
	 \}$, and $\inj(M)$ is the injectivity radius of $M$.  Then the 
	 result follows from (\cite{Chavel}, p.119,123). \\
	 \noindent $(ii)$ 
	 \begin{align*}
		 \lim_{r\to0} \frac{\b_2(r)}{\b_1(r)} & = \lim_{r\to0} 
		 \frac{V_{c_2}(r)}{V_{c_1}(r)} = \lim_{r\to0} \frac{ \int_0^r 
		 S_{c_2}(t)^{n-1}dt }{ \int_0^r S_{c_1}(t)^{n-1}dt } \\
		 & = \left( \lim_{r\to0} \frac{S_{c_2}(r)}{S_{c_1}(r)} 
		 \right)^{n-1} = 1
	 \end{align*}
	 where (cfr.  \cite{Chavel}, formulas (2.48), (3.24), (3.25)) 
	 $V_\d(r) = \frac{n\sqrt{\pi} }{ \G(n/2+1)} \int_0^r 
	 S_\d(t)^{n-1}dt$, and
	 $$
	 S_\d(r) := 
	 \begin{cases}
		 \frac1{\sqrt{-\d}}\sinh(r\sqrt{-\d}) & \d<0 \\
		 r& \d=0 \\
		 \frac1{\sqrt{\d}}\sin(r\sqrt{\d}) &  \d>0. 
	 \end{cases}
	 $$
 \end{proof}
 
 \medskip 

 Let $M$ be a complete Riemannian manifold, and recall (\cite{Wolf}) 
 that $\D_{p}:=(d+d^{*})^{2}|_{L^{2}(\La^{p}T^{*}M)}$, the $p$-th 
 Laplacian on $M$, is essentially self-adjoint and positive, and the 
 semigroup $\e{-t\D_{p}}$ has a C$^\infty$ kernel, $H_{p}(t,x,y)$, on 
 $(0,\infty)\times M\times M$, called the $p$-th heat kernel.  Let us 
 mention the following result, which will be useful in the sequel.

 \begin{Prop}\label{2.1.5} 
	 Let $M$ be an $n$-dimensional complete Riemannian manifold with 
	 C$^{\infty}$ bounded geometry, then for all $T>0$, there are $c, 
	 c'>0$, s.t., for $0<t\leq T$,
	 \begin{align*}
		 |H_{p}(t,x,y)| & \leq c\ t^{-n/2-1} 
		 \exp \left(\frac{-c'\d(x,y)^{2}}{t}\right)  \\
		 |\nabla_{x}H_{p}(t,x,y)| & \leq c\ t^{-n/2-3/2} 
		 \exp \left(\frac{-c'\d(x,y)^{2}}{t}\right) 
	 \end{align*}
	 where we denoted with $\d$ the metric induced on $M$ by $g$.  As 
	 a consequence $H_{p}(t,\cdot,\cdot)$ is uniformly continuous on a 
	 neighborhood of the diagonal of $M\times M$.
 \end{Prop}
 \begin{proof}
	 The estimates are proved in \cite{BE}.  For the last statement, 
	 for any $\d_{0}<\min\{1,\inj(M),\frac{\pi}{ \sqrt{c_{1}} } \}$, 
	 $x\in M$, $y\in B(x,\d_{0})$, we have 
	 $|H_{p}(t,x,y)-H_{p}(t,x,x)|\leq \sup|\nabla_{y}H_{p}(t,x,y)| 
	 \d(x,y)$, and we get the uniform continuity.
 \end{proof}
 
 \medskip

 \subsection{The C$^*$-algebra of almost local 
 operators}\label{subsec:almostlocal}

 Let $F$ be a finite dimensional Hermitian vector bundle over $M$, and 
 let $L^{2}(F)$ be the Hilbert space completion of the smooth sections 
 with compact support of $F$ w.r.t. the scalar product $\langle 
 s_{1},s_{2} \rangle := \int_{M} \langle s_{1x},s_{2x} \rangle 
 dvol(x)$.
 
 Recall \cite{RoeCoarse} that an operator $A\in\cb(L^2(F))$ has finite 
 propagation if there is a constant $u_A>0$ s.t. for any compact 
 subset $K$ of $M$, any $\f\in L^2(F)$, $\supp\f\subset K$, we have 
 $\supp A\f \subset Pen^+(K,u_A) := \{ x\in M : \d(x,K)\leq u_A \}$.
 
 Let us denote by $\ca_0\equiv \ca_{0}(F)$ the set of finite 
 propagation operators.  $\ca_0$ may be characterized as follows

 \begin{Prop}\label{3.1.1} 
	 \itm{i} $A\in\ca_0$ iff, for any measurable set $\O$,
	 $AE_\O=E_{Pen^+(\O,u_A)}AE_\O$, where $E_X$ is the multiplication
	 operator by the characteristic function of the set $X$; 
	 \itm{ii} $A\in\ca_0$ iff, for any functions $\f,\ \psi\in L^2(F)$
	 with $\d(\supp\psi, \supp\f)> u_A$, one has $(\f,A\psi)=0$.
 \end{Prop}
 \begin{proof} Properties $(i)$ and $(ii)\ (\imply)$ are obvious. \\
	 $(ii)\ (\coimply)$ The hypothesis implies that $\supp A\psi 
	 \subset M\setminus\supp\f$ for all $\f$ s.t. $\supp\f 
	 \subset M\setminus Pen^+(\supp\psi,u_A)$.  The thesis follows.  
 \end{proof}

 \begin{Prop}\label{3.1.2} 
	 The set $\ca_0(F)$ of finite propagation operators is a 
	 $^*$-al\-ge\-bra with identity.
 \end{Prop}
 \begin{proof}
	 Let $K$ be a compact subset of $M$, $\f\in L^2(F)$, 
	 $\supp\f\subset K$, and $A,\ B\in \ca_0$.  Then $\supp(A+B)\f 
	 \subset \supp A\f \cup \supp B\f$, hence $u_{A+B}= u_A\vee u_B$ 
	 is the requested constant.  Moreover $\supp(AB)\f \subset 
	 Pen^+(\supp B\f,u_A) \subset Pen^+(K, u_A + u_B)$, so that we may 
	 set $u_{AB} =u_A+u_B$.  \\
	 As $(A^*\psi,\f) = (\psi,A\f)=0$ for all $\f,\ \psi\in L^2(F)$, 
	 with $\d(\supp\psi, \supp\f) > u_A$, that is $\supp \f \cap 
	 Pen^+(\supp \psi,u_A) = \emptyset$, we get $\supp A^*\psi \subset 
	 Pen^+(\supp\psi,u_A)$, which implies $u_{A^*}\leq u_A$, and 
	 exchanging the roles of $A,\ A^*$, we get $u_A= u_{A^*}$.  
 \end{proof}
 
 \medskip
 
 The norm closure of $\ca_0$ will be denoted by $\ca\equiv \ca(F)$ and 
 will be called the C$^*$-algebra of almost local operators on 
 $L^{2}(F)$.  Now we show that Gaussian decay for the kernel of a 
 positive operator $A$ is a sufficient condition for $A$ to belong to 
 $\ca$.

 \begin{Thm}\label{3.1.3} 
	 Let $M$ be a complete Riemannian $n$-manifold of C$^{\infty}$ 
	 bounded geometry.  If $A$ is a bounded self-adjoint operator on 
	 $L^2(F)$, with kernel $a(x,y)\in L(F_{y},F_{x})$, and there are 
	 positive constants $c,\ \a,\ \d_0$ s.t., for $\d(x,y)\geq\d_0$, 
	 $a(x,y)$ is measurable and
	 $$ 
	 |a(x,y)|\leq c\ \e{-\a\d(x,y)^2}
	 $$
	 then $A\in\ca$.
 \end{Thm}

 In order to prove the theorem, we need some lemmas.
 
 \begin{Lemma}\label{3.1.4} 
	 Let $A$ be a bounded self-adjoint operator on $L^2(F)$, with 
	 measurable kernel.  Then
	 $$
	 \|A\| \leq \sup_{x\in M} \int_M |a(x,y)|dy
	 $$
 \end{Lemma}
 \begin{proof}
	 Since $A$ is self-adjoint, $a(x,y)$ is symmetric, hence
	 \begin{align*}
		 \|A\|_{1\to1} & = \sup \{|(f,Ag)| : f\in L^\infty(F),\ 
		 \|f\|_\infty=1,\ g\in L^1(F),\ \|g\|_1=1 \} \\
		 & \leq \sup_{x\in M} \int_M |a(y,x)|dy = 
		 \|A\|_{\infty\to\infty}
	 \end{align*}
	 The thesis easily follows from Riesz-Thorin interpolation 
	 theorem.
 \end{proof}
 
 \begin{Lemma}\label{3.1.5} 
	 Let $\f:[0,\infty)\to[0,\infty)$ be a non-increasing measurable 
	 function.  Then, using notation of subsection 1.1,
	 $$
	 \sup_{x\in M} \int_M \f(\d(x,y))dy \leq C_n\ \int_0^\infty 
	 \f(r)S_{c_2}(r)^{n-1}dr
	 $$
	 where $C_n:= \frac{n\sqrt{\pi}}{\G(n/2+1)}$, and $S_{c_2}(r):= 
	 \frac1{\sqrt{-c_2}}\sinh(r\sqrt{-c_2})$.
 \end{Lemma} 
 \begin{proof}
	 From Theorem \ref{2.1.2} we get $V(x,r)\leq C_n\ \int_0^r 
	 S_{c_2}(t)^{n-1}dt$.  Then
	 \begin{align*}
		 \int_M \f(\d(x,y))dy & = \int_0^\infty \f(r)dV(x,r) \\
		 & \leq C_n\ \int_0^\infty \f(r) S_{c_2}(r)^{n-1}dr
	 \end{align*}
	 where the equality is in $e.g.$ (\cite{HeSt}, Theorem 12.46), and 
	 the inequality holds because $\f$ is non-increasing and positive, 
	 and $V(x,0)=0$.  
 \end{proof}
 
 \medskip

 \noindent {\it Proof of Theorem \ref{3.1.3}.} 
 Let $\r>\d_0$, and decompose $A=A_\r+A'_\r$, with $a_\r(x,y):= 
 a(x,y)\chi_{[0,\r]}(\d(x,y))$.  Then $A_\r\in\ca_0$, and 
 $|a'_\r(x,y)|\leq c'\f(\d(x,y))$, where
 $$
 \f(r):=
 \begin{cases}
	 \e{-\a\r^2} &  0\leq r<\r \\
	 \e{-\a r^2} &  r\geq\r. 
 \end{cases}
 $$
 By Lemmas \ref{3.1.4}, \ref{3.1.5} we get 
 \begin{align*}
	 \|A-A_\r\| = \|A'_\r\| & \leq \sup_{x\in M} \int_M |a'_\r(x,y)|dy 
	 \\
	 & \leq c'\ \sup_{x\in M} \int_M \f(\d(x,y))dy \\
	 & \leq c' \int_0^\infty \f(r)S_{c_2}(r)^{n-1}dr \\
	 & \leq c'\ \e{-\a\r^2} \int_0^\r S_{c_2}(r)^{n-1}dr + c'' 
	 \int_\r^\infty \e{-\a r^2+(n-1)r\sqrt{-c_2}}dr\\
	 & \to 0,\quad \r\to\infty
 \end{align*}
 and the thesis follows.
 \qed
 
 \medskip
 
 Applying the previous Theorem we conclude that $C_0$ functional 
 calculus of the Laplace operator belongs to $\ca$.

 \begin{Cor}\label{3.1.6} 
	 Let $M$ be a complete Riemannian manifold of C$^{\infty}$ bounded 
	 geometry.  Then $\f(\D_{p})\in\ca(\La^{p}T^{*}M)$, for any $\f\in 
	 C_0([0,\infty))$.
 \end{Cor}
 \begin{proof}
	 By Proposition \ref{2.1.5} and Theorem \ref{3.1.3} we obtain that 
	 $\e{-t\D_{p}}\in\ca$, for any $t>0$.  Since $\{ \e{-t\l} 
	 \}_{t>0}$ generates a $^*$-algebra of $C_0([0,\infty))$ which 
	 separates points, the thesis follows by Stone-Weierstrass 
	 theorem.
 \end{proof}

 \begin{Prop}
    \itm{i} $\ca$ contains all compact operators, so that $\ca'' = 
\cb(\ch)$
    \itm{ii} $\ca=\cb(\ch)$ iff $M$ is a closed manifold.
 \end{Prop}
 \begin{proof}
	 $(i)$ It suffices to prove that $\ca$ contains all 
	 one-dimensional projections.  To begin with, let $e$ be the 
	 projection operator onto the multiples of $\f\in 
	 C_{c}(\La^{p}T^{*}M)$ with L$^{2}$-norm $1$, so that its Schwartz 
	 kernel is $K(x,y):=\f(x)\otimes\ov{\f(y)}$.  Consider the 
	 operators $e_{n}$ with Schwartz kernel $K_{n}(x,y) := 
	 \chi_{B(o,n)}(x)\f(x)\otimes\chi_{B(o,n)}(y)\ov{\f(y)}$.  Then it 
	 follows from Lemma \ref{3.1.4} that
	 \begin{align*}
		\|e-e_{n}\| & \leq \sup_{x\in M} \int_{M} |\f(x)\f(y) - 
		\chi_{B(o,n)}(x)\f(x)\chi_{B(o,n)}(y)\f(y)| dy \\
		& = \sup_{x\in M} |\f(x)| \int_{M} |\f(y)|(1 - 
		\chi_{B(o,n)}(x)\chi_{B(o,n)}(y)) dy \\
		& \leq \max\{ \sup_{x\in B(o,n)} |\f(x)| \int_{B(o,n)^{c}} |\f(y)| 
		dy, \sup_{x\in B(o,n)^{c}} |\f(x)| \int_{M} |\f(y)| dy \} \to 0.
	 \end{align*}
	 Therefore $e\in\ca$.  Let now $\f\in L^{2}(\La^{p}T^{*}M)$, with 
	 $\|\f\|_{2}=1$, and let $\{\f_{n}\}\subset C_{c}(\La^{p}T^{*}M)$ 
	 s.t. $\|\f-\f_{n}\|_{2}<\frac{1}{n}$.  Denote by $e$, resp.  
	 $e_{n}$, the projection operator onto the multiples of $\f$, 
	 resp.  $\f_{n}$.  Then $e_{n}\to e$, so that $e\in\ca$. 
	 \\
	 $(ii)$ If $M$ is closed the statement is trivial. Conversely
	 let $M$ be open and $F$ be the trivial bundle, as we may
	 reduce to this case.  Then choose a sequence $\{x_n\}$ in $M$
	 s.t. $\d(x_i,x_0)\geq 3i$, and $\d(x_i,x_j)\geq 3$, $i\neq
	 j$. Then $f_i := \frac {\chi_{B(x_i,1)}} {V(x_i,1)}$ is an
	 orthonormal set in $L^2(M)$. The operator $T:=(f_0,\cdot)
	 \sum_{i=1}^\infty f_i$ has norm 1, but $T\not\in\ca$, because
	 $\|T-A\|\geq 1$ for any $A\in\ca_0$.
 \end{proof} 

 \subsection{A functional described by J.~Roe}
 \label{subsec:roe}

 In the rest of this paper $M$ is a complete Riemannian $n$-manifold
 of C$^\infty$ bounded geometry as in Definition \ref{2.1.1}, that we
 assume endowed with a regular exhaustion $\ck$ \cite{Roe1}, that is
 an increasing sequence $\{K_{n}\}$ of compact subsets of $M$, whose
 union is $M$, and s.t., for any $r>0$
 $$									
 \lim_{n\to\infty} \frac{vol(K_{n}(r)) }{vol(K_{n}(-r))} =1,
 $$
 where we set, here and in the following, $K(r)\equiv Pen^{+}(K,r):=
 \{x\in M: \d(x,K)\leq r\}$, and $K(-r)\equiv Pen^{-}(K,r):=$ the
 closure of $M\setminus Pen^{+}(M\setminus K,r)$.  Observe that, as
 $M$ is complete, $Pen^+(K,r)$ coincides with the closure of $\{ x\in
 M : \d(x,K) < r \}$, which is the original definition of Roe.
  
 \begin{Lemma}\label{l:penombra} 
	 Let $K$ be a compact subset of $M$, then
	 \itm{i} $K(-r_{2})\subset K\subset K(r_{1})$, for any $r_1,r_2>0$
	 \itm{ii} $\{x\in M: \d(x,M\setminus K)<r\} \subset$ Interior of 
 $Pen^+(M\setminus K,r) \equiv M\setminus K(-r)$ 
	 \itm{iii} $Pen^{+}(K(r_{1})\setminus K(-r_{2}),R) \subset 
	 K(r_{1}+R+\eps)\setminus K(-r_{2}-R-\eps)$, for any 
	 $r_1,r_2,R,\eps>0$.
 \end{Lemma}
 \begin{proof}
	 $(ii)$ If $\d(x,M\setminus K)<r$, there is $z\in M\setminus K$ 
	 s.t. $\d(x,z)<r$, so that $x$ belongs to the interior of 
	 $Pen^{+}(M\setminus K,r)$, which is the complement of $K(-r)$.\\
	 $(iii)$ Indeed if $x\in Pen^{+}(K(r_{1})\setminus K(-r_{2}),R)$, 
	 then for any $\eps>0$ there is $x_{\eps}\in K(r_{1})\setminus 
	 K(-r_{2})$ with $\d(x,x_\eps)<R+\eps/2$.  Therefore, on the one 
	 hand, $\d(x,K) \leq R + \eps/2 + r_{1}$, which implies $x\in 
	 K(r_{1}+R+\eps)$.  On the other hand, as $x_{\eps}\not\in 
	 K(-r_{2})$, there is $y_{\eps}\in M\setminus K$ s.t. 
	 $\d(y_{\eps},x_{\eps})<r_{2}+\frac{\eps}2$, hence 
	 $\d(x,y_{\eps})\leq \frac{\eps}2 + R + r_{2}+\frac{\eps}2$ and 
	 $x\not\in K(-r_{2}-R-\eps)$.
 \end{proof}
  
 Following Moore-Schochet \cite{Moore-Schochet}, we recall that an 
 operator $T$ on $L^2(F)$ is called locally trace class if, for any 
 compact set $K\subset M$, $E_KTE_K$ is trace class, where $E_K$ 
 denotes the projection given by the characteristic function of $K$.  
 It is known that the functional $\m_T(K):= Tr(E_KTE_K)$ extends to a 
 Radon measure on $M$.  To state next definition we need some 
 preliminary notions.

%%%fin qui

 \begin{Dfn} 
	 Define $\cj_{0+}\equiv \cj_{0+}(F)$ as the set of positive 
	 locally trace class operators $T$, such that 
	 \itm{i} there is $c>0$ s.t. $\m_T(K_{n}) \leq c\ vol(K_{n})$, 
	 asymptotically, 
	 \itm{ii} $\lim_{n\to\infty}\frac{\m_{T}(K_{n}(r_{1})\setminus 
	 K_{n}(-r_{2}))}{vol(K_{n})} = 0$.
 \end{Dfn}
 
 \begin{Lemma}\label{3.2.3} 
	 $\cj_{0+}$ is a hereditary (positive) cone in $\cb(L^2(F))$.
 \end{Lemma}
 \begin{proof} Linearity follows by $\mu_{A+B}=\mu_A+\mu_B$.
	 If $T\in\cj_{0+}$, and $0\leq A\leq T$, then $Tr(BAB^*)\leq 
	 Tr(BTB^*)$, for any $B\in\cb(L^2(F))$, and the thesis follows. 	 
 \end{proof}
 
 \begin{rem}\label{rem:laplacian} 
	 The hereditary cone $\cj_{0+}$ depends on the exhaustion $\ck$, 
	 however it contains a (hereditary) subcone, given by the 
	 operators $T$ for which there is $c>0$ such that 
	 $\m_T(\Omega)\leq c\ vol(\Omega)$ for any measurable set 
	 $\Omega$.  Proposition~\ref{2.1.5} implies that the operator 
	 $e^{-t\Delta_p}$ belongs to the subcone, hence to 
	 $\cj_{0+}(\La^{p}T^*{M})$.
 \end{rem}

 Recall \cite{Roe1} that $\cu_{-\infty}(F)$ is the set of uniform 
 operators of order $-\infty$.
 
 \begin{Prop}\label{prop:cu}
 	$\cu_{-\infty}(F)_{+} \subset \cj_{0+}(F)$.
 \end{Prop}
 \begin{proof}
 	Let $A\in \cu_{-\infty}(F)$, so that $Au(x) = \int_{M} a(x,y)u(y) 
 	dy$, with $a\in C^{\infty}(F\otimes F)$ is a smoothing kernel, and 
 	is uniformly bounded together with all its covariant derivatives 
 	(\cite{Roe1}, 2.9). Then for any Borel set $\O\subset M$, 
 	$\m_{A}(\O) = Tr(E_{\O}AE_{\O}) = \int_{\O} tr(a(x,x))dx \leq c\ 
 	vol(\O)$, and the result easily follows.
 \end{proof}
 
 If $\om$ is a state on $\ell^\infty(\bn)$ vanishing on infinitesimal 
 sequences, we use in the following the notation $\lo a_{n}:= 
 \om(\{a_{n}\})$, for any $\{a_{n}\}\in\ell^\infty(\bn)$.  Consider the 
 weight $\f\equiv \f_{\ck,\om}$ on $\cb(L^2(F))_+$ given by
 $$
 \f(A):= 
 \begin{cases}
	 \lo \frac{\m_A(K_{n})}{ vol(K_{n})}  &  A\in\cj_{0+} \\
	 +\infty &  A\in\cb(L^2(F))_+\setminus\cj_{0+}.
 \end{cases}
 $$
 Observe that the functional $\f$ is the functional defined by Roe 
 in \cite{Roe1}, but for the domain.
  
 \begin{Prop}
	For any $A\in\cu_{-\infty}(F)_+$, $\f(A) = \lo \frac{\int_{K_{n}} 
	tr(a(x,x)) dx}{vol(K_{n})}$, which is Roe's definition in {\rm 
	\cite{Roe1}}.
 \end{Prop}
 \begin{proof}
 	Follows easily from the proof of Proposition \ref{prop:cu}.
 \end{proof}
 
 \begin{Lemma}\label{l:phiomega}
	 If $A\in\cj_{0+}$ then 
	 \begin{equation*}
		\f(A) = \lo \frac{\m_A(K_{n}(r_{1}))}{ vol(K_{n}(r_{2}))}
	 \end{equation*}
	 for any $r_1,r_2\in\br$.
 \end{Lemma}
 \begin{proof}
	 Indeed, if $r_{1}\geq0$, we get
	 \begin{equation*}
		 \lo\frac{\m_A(K_{n}(r_{1}))}{ vol(K_{n}(r_{2}))} =\lo\left( 
		 \frac{\m_A(K_{n})}{ vol(K_{n})} + 
		 \frac{\m_A(K_{n}(r_{1})\setminus K_{n})}{ vol(K_{n})}\right) 
		 \frac{vol(K_{n})}{ vol(K_{n}(r_{2}))} = \f(A)
	 \end{equation*}
	 whereas, if $r_{1}<0$, we get
	 \begin{equation*}
		\lo\frac{\m_A(K_{n}(r_{1}))}{ vol(K_{n}(r_{2}))} =\lo\left( 
		\frac{\m_A(K_{n})}{ vol(K_{n})} - \frac{\m_A(K_{n}\setminus 
		K_{n}(r_{1}))}{ vol(K_{n})}\right) \frac{vol(K_{n})}{ 
		vol(K_{n}(r_{2}))} =\f(A).
	 \end{equation*}
 \end{proof}
 
 The algebra $\ca$, being a C$^*$-algebra, contains many unitary 
 operators, and is indeed generated by them.  The algebra $\ca_0$ may 
 not, but all unitaries in $\ca$ may be approximated by elements in 
 $\ca_0$.  Such approximants are $\d$-unitaries, according to the 
 following
 
 \begin{Dfn}\label{3.2.5} 
	 An operator $U\in\cb(L^2(F))$ is called $\d$-unitary, $\d>0$, if 
	 $\|U^*U-1\|<\d$, and $\|UU^*-1\|<\d$.
 \end{Dfn}

 Let us denote with $\cu_\d$ the set of $\d$-unitaries in $\ca_0$ and 
 observe that, if $\d<1$, $\cu_\d$ consists of invertible operators, 
 and $U\in\cu_\d$ implies $U^{-1}\in\cu_{\d/(1-\d)}$.

 \begin{Prop}\label{3.2.6} 
	 The weight $\f$ is $\eps$-invariant for $\d$-unitaries in 
	 $\ca_0$, namely, for any $\eps\in(0,1)$, there is $\d>0$ s.t., 
	 for any $U\in\cu_\d$, and $A\in\ca_+$,
	 $$
	 (1-\eps)\f(A) \leq \f(UAU^*) \leq (1+\eps)\f(A) .
	 $$
 \end{Prop}

 \begin{Lemma}\label{3.2.7} 
	 If $T\in\cj_{0+}$, then $ATA^*\in\cj_{0+}$ for all $A\in\ca_0$.
 \end{Lemma}
 \begin{proof}
	 First observe that for any Borel set $\Omega\subset M$ we have
	 \begin{align*}
		 \m_{ATA^*}(\Omega)
		 & = Tr(E_\Omega ATA^*E_\Omega) \\
		 & = Tr(E_\Omega AE_{\Omega(u_A)}TE_{\Omega(u_A)}A^*E_\Omega) \\  
		 & \leq \|A^*E_\Omega A\| Tr(E_{\Omega(u_A)}TE_{\Omega(u_A)}) \\  
		 & \leq \|A\|^2 \m_{T}(\Omega(u_A))
	 \end{align*}
	 so that 
	 \begin{align*}
		 \frac{\m_{ATA^*}(K_{n})}{vol(K_{n})}
		 & \leq \|A\|^2 \frac{\m_{T}(K_{n}(u_A))}{vol(K_{n})} \\
		 & = \|A\|^2 \frac{\m_{T}(K_{n})}{vol(K_{n})} + 
		 \|A\|^2\frac{\m_{T}(K_{n}(u_A)\setminus K_{n})}{vol(K_{n})}
	 \end{align*}
	 which is asymptotically bounded.  Now observe that, by Lemma 
	 \ref{l:penombra} $(iii)$, it follows
	 \begin{align*}
		 \frac{\m_{ATA^*}(K_{n}(r_{1})\setminus K_{n}(-r_{2}))}{vol 
		 K_{n}} & \leq \|A\|^{2} 
		 \frac{\m_{T}(Pen^{+}(K_{n}(r_{1})\setminus 
		 K_{n}(-r_{2}),u_A))}{vol K_{n}} \\
		 & \leq \|A\|^{2} \frac{\m_{T}(K_{n}(r_{1}+u_A+\eps)\setminus 
		 K_{n}(-r_{2}-u_A-\eps))}{vol K_{n}} \to 0,
	 \end{align*}
	 the thesis follows.
 \end{proof}
 
 \medskip
  
 \noindent {\it Proof of Proposition \ref{3.2.6}.}
 Assume $A\in \cj_{0+}\cap\ca_{+}$, then $UAU^{*}\in\cj_{0+}$ and, by 
 Lemma \ref{l:phiomega},
 \begin{align*}
	 \f(UAU^*) & = \lo \frac{\m_{UAU^{*}}(K_{n})}{ vol(K_{n})} \\
	 & \leq \|U\|^2 \lo\left(\frac{\m_A(K_{n}(u_{U})) }{ vol(K_{n})} 
	 \right)\\
	 & \leq (1+\d) \f(A).
 \end{align*}
 Choose now $\d<\eps/2$, and $U\in\cu_\d$, so that 
 $U^{-1}\in\cu_{2\d}$, and $\f(UAU^*)\leq (1+\d) \f(A)< (1+\eps) 
 \f(A)$.  Replacing $A$ with $UAU^*$, and $U$ with $U^{-1}$, we obtain
 $$
 \f(A)\leq \|U^{-1}\|^2 \f(UAU^*)\leq (1+2\d)\f(UAU^*) < 
 (1+\eps)\f(UAU^*)
 $$
 and the thesis easily follows. \\
 Assume now $A\in\ca_{+}\setminus\cj_{0+}$, then 
 $\f(A)=+\infty=\f(UAU^{*})$, because otherwise $UAU^{*}\in\cj_{0+}$, 
 so that $A=U^{-1}(UAU^{*})(U^{-1})^{*}\in\cj_{0+}$, which is absurd.  
 \qed
 
 \medskip
 
 Finally we observe that, from the proof of Lemma \ref{3.2.7} the 
 following is immediately obtained
 
 \begin{Lemma}\label{3.2.8} 
	 If $A\in\ca_0$ and $\|A\|\leq 1$, then $\f(ATA^*)\leq \f(T)$, for 
	 any $T\in\cj_{0+}$.
 \end{Lemma}

 \subsection{A construction of semicontinuous traces on 
 C$^*$-algebras} \label{subsec:semicontinuous}

 The purpose of this subsection is to show that the 
 lower-semicontinuous semifinite regularisation of the functional 
 $\f|_\ca$ of the previous subsection gives a trace, namely a 
 unitarily invariant weight on $\ca$.  It turns out that this 
 procedure can be applied to any weight $\t_0$, on a unital 
 C$^*$-algebra $\ca$, which is $\eps$-invariant for $\d$-unitaries of 
 a dense $^*$-subalgebra $\ca_0$.  The particular case of the 
 functional $\f|_\ca$ is treated in the next subsection.  First we 
 observe that, with each weight on $\ca$, namely a functional 
 $\t_0:\ca_+\to[0,\infty]$, satisfying the property $\t_0(\l 
 A+B)=\l\t_0(A)+\t_0(B)$, $\l>0$, $A,\ B\in\ca_+$, we may associate a 
 (lower-)semicontinuous weight $\t$ with the following procedure
 \begin{equation}\label{e:weight}
	 \t(A):= \sup \{ \psi(A): \psi\in\ca^*_+,\ \psi\leq\t_0 \}
 \end{equation}
 Indeed, it is known that \cite{Combes,Stratila}
 $$
 \t(A) \equiv \sup_{\psi\in\cf(\t_0)} \psi(A) 
 $$
 where $\cf(\t_0):= \{ \psi\in\ca^*_+ : \exists\ \eps>0,\ 
 (1+\eps)\psi<\t_0 \}$.  Moreover the following holds

 \begin{Thm}\label{3.3.1} {\rm \cite{QV}} 
	 The set $\cf(\t_0)$ is directed, namely, for any $\psi_1,\ 
	 \psi_2\in\cf(\t_0)$, there is $\psi\in\cf(\t_0)$, s.t. $\psi_1,\ 
	 \psi_2\leq \psi$.
 \end{Thm}

 From this theorem easily follows

 \begin{Cor}\label{3.3.2} 
	 Let $\t_0$ be a weight on the C$^*$-algebra $\ca$, and $\t$ be 
	 defined as in $(\ref{e:weight})$.  Then 
	 \itm{i} $\t$ is a semicontinuous weight on $\ca$ 
	 \itm{ii} $\t=\t_0$ iff $\t_0$ is semicontinuous.  
	 \itm{iii} The domain of $\t$ contains the domain of $\t_0$.
	 \\
	 The weight $\t$ will be called the semicontinuous regularization 
	 of $\t_0$.
 \end{Cor}
 \begin{proof}
	 $(i)$ From Theorem \ref{3.3.1}, $\t(A) = \sup_{\psi\in\cf(\t_0)} 
	 \psi(A) = \lim_{\psi\in\cf(\t_0)} \psi(A)$, whence linearity and 
	 semicontinuity of $\t$ easily follow.  \\
	 $(ii)$ is a well known result by Combes \cite{Combes}.\\
	 $(iii)$ Immediately follows from the definition of $\t$.
 \end{proof}

 \begin{Prop}\label{3.3.3} 
	 Let $\t_0$ be a weight on $\ca$ which is $\eps$-invariant by 
	 $\d$-unitaries in $\ca_0$ (as in Proposition \ref{3.2.6}).  Then 
	 the associated semicontinuous weight $\t$ satisfies the same 
	 property.
 \end{Prop}
 \begin{proof}
	 Fix $\eps<1$ and choose $\d\in(0,1/2)$, s.t. $U\in\cu_\d$ implies 
	 $|\t_0(UAU^*)-\t_0(A)|<\eps \t_0(A)$, $A\in\ca_+$.  Then, for any 
	 $U\in\cu_{\d/2}$ and any 
	 $\psi\in\ca^*_+$, $\psi\leq\t_0$, we get
	 $$
	 \psi\circ adU(A)\leq
	 \t_0(UAU^*)\leq (1+\eps)\t_0(A),
	 $$
	 for $A\in\ca_+$, $i.e.$ $(1+\eps)^{-1}\psi\circ adU\leq \t_0$.  
	 Then
	 \begin{align*}
		 \t(UAU^*) 
		 & = (1+\eps) \sup_{\psi\leq\t_0} (1+\eps)^{-1}\psi\circ adU(A) \\
		 & \leq (1+\eps) \sup_{\psi\leq\t_0} \psi(A) \\
		 & = (1+\eps)\t(A).
	 \end{align*}
	 Since $U^{-1}\in\cu_\d$, replacing $U$ with $U^{-1}$ and $A$ with 
	 $UAU^*$, we get $\t(A)\leq (1+\eps)\t(UAU^*)$.  Combining the 
	 last two inequalities, we get the result.  
 \end{proof}

 \begin{Prop}\label{3.3.4} 
	 The semicontinuous weight $\t$ is a trace on $\ca$, namely, 
	 setting $\cj_+:= \{A\in\ca_+: \t(A)<\infty\}$, and extending $\t$ 
	 to the linear span $\cj$ of $\cj_+$, we get
	 \itm{i} $\cj$ is an ideal in $\ca$ 
	 \itm{ii} $\t(AB)=\t(BA)$, for all $A\in\cj$, $B\in\ca$.
 \end{Prop}
 \begin{proof}
	 $(i)$ Let us prove that $\cj_+$ is a unitary invariant face in 
	 $\ca_+$, and it suffices to prove that $A\in\cj_+$ implies 
	 $UAU^*\in\cj_+$, for all $U\in\cu(\ca)$, the set of unitaries in 
	 $\ca$.  Suppose on the contrary that there is $U\in\cu(\ca)$ s.t. 
	 $\t(UAU^*)=\infty$.  Then there is $\psi\in\ca^*_+$, 
	 $\psi\leq\t_0$, s.t. $\psi(UAU^*) > 2\t(A)+2$.  Then we choose 
	 $\d<3$ s.t. $V\in\cu_\d$ implies $\t(VAV^*)\leq 2\t(A)$, and an 
	 operator $U_0\in\ca_0$ s.t. $\|U-U_0\|< \min\{ \frac{\d}{3}, 
	 \frac1{3\|A\|\|\psi\|} \}$.  The inequalities
	 $$
	 \|U_0U_0^*-1\| = \|U^*U_0U_0^*-U^*\| \leq 
	 \|U^*U_0-1\|\|U_0^*\|+\|U_0^*-U^*\| < \d
	 $$
	 and analogously for $\|U_0^*U_0-1\|<\d$, show that 
	 $U_0\in\cu_\d$.  Then, since $|\psi(U_0AU_0^*)-\psi(UAU^*)|\leq 
	 3\|\psi\|\|A\|\|U-U_0\|<1$, we get
	 $$
	 2\t(A)\geq \t(U_0AU_0^*) \geq \psi(U_0AU_0^*) \geq
	 \psi(UAU^*) - 1 \geq 2\t(A)+1
	 $$
	 which is absurd. \\
	 $(ii)$ We only have to show that $\t$ is unitary invariant.  Take 
	 $A\in\cj_+$, $U\in\cu(\ca)$.  For any $\eps>0$ we may find a 
	 $\psi\in\ca^*_+$, $\psi\leq\t_0$, s.t. 
	 $\psi(UAU^*)>\t(UAU^*)-\eps$, as, by $(i)$, $\t(UAU^*)$ is 
	 finite.  Then, arguing as in the proof of $(i)$, we may find 
	 $U_0\in\ca_0$, so close to $U$ that
	 \begin{align*}
		 & |\psi(U_0AU_0^*)-\psi(UAU^*)|<\eps \\
		 & (1-\eps)\t(A)\leq \t(U_0AU_0^*) \leq (1+\eps)\t(A). 
	 \end{align*}
	 Then
	 \begin{align*}
		 \t(A) & \geq \frac1{1+\eps}\ \t(U_0AU_0^*) \geq 
		 \frac1{1+\eps}\ \psi(U_0AU_0^*) \\
		 & \geq \frac1{1+\eps}\ (\psi(UAU^*) -\eps) \geq 
		 \frac1{1+\eps}\ (\t(UAU^*) -2\eps).
	 \end{align*}
	 By the arbitrariness of $\eps$ we get $\t(A)\geq \t(UAU^*)$.  
	 Replacing $A$ with $UAU^*$, we get the thesis.  
 \end{proof}
  
 The second regularization we need turns $\t$ into a (lower 
 semicontinuous) semifinite trace, namely guarantees that
 $$
 \t(A) = \sup\{ \t(B) : 0\leq B\leq A,\ B\in\cj_+ \}
 $$
 for all $A\in\ca_+$.  In particular the semifinite regularization 
 coincides with the original trace on the domain of the latter.  This 
 regularization is well known (see $e.g.$ \cite{DixmierC}, Section 6), 
 and amounts to represent $\ca$ via the GNS representation $\pi$ 
 induced by $\t$, define a normal semifinite faithful trace $tr$ on 
 $\pi(\ca)''$, and finally pull it back on $\ca$, that is 
 $tr\circ\pi$.  It turns out that $tr\circ\pi$ is (lower 
 semicontinuous and) semifinite on $\ca$, $tr\circ\pi\leq\t$, and 
 $tr\circ\pi(A)=\t(A)$ for all $A\in\cj_+$, that is $tr\circ\pi$ is a 
 semifinite extension of $\t$, and $tr\circ\pi=\t$ iff $\t$ is 
 semifinite.

 We still denote by $\t$ its semifinite extension.  As follows from 
 the construction, semicontinuous semifinite traces are exactly those 
 of the form $tr\circ\pi$, where $\pi$ is a tracial representation, 
 and $tr$ is a n.s.f. trace on $\pi(\ca)''$.

 \subsection{The regularized trace on the C$^*$-algebra of almost 
 local operators} \label{subsec:regularized}

 Now we apply the regularization procedure described in the previous 
 subsection to Roe's functional.  Let us remark that the 
 semicontinuous regularization of the weight $\f|_\ca$ is a trace in 
 the sense of property $(ii)$ of Proposition~\ref{3.3.4}, which is 
 stronger then the trace property in \cite{Roe1}.  First we observe 
 that $\f|_\ca$ is not semicontinuous.

 \begin{Prop}\label{3.4.1} 
	 The set $\cn_0:=\{T\in\ca_+: \f(T)=0\}$ is not closed.  In 
	 particular, there are operators $T\in\ca_+$ s.t. $\f(T)=1$ but 
	 $\t(T)=0$ for any (lower-)semicontinuous trace $\t$ dominated by 
	 $\f|_\ca$.
 \end{Prop}
 \begin{proof}
	 Recall from Lemma \ref{2.1.2}$(i)$ that there are positive real 
	 functions $\b_1,\ \b_2$ s.t. $0<\b_1(r)\leq V(x,r) \leq \b_2(r)$, 
	 for all $x\in M$, $r>0$, and $\lim_{r\to0}\b_2(r)=0$.  Therefore 
	 we can find a sequence $r_n\dec0$ s.t. $\sum_{n=1}^\infty 
	 \b_2(r_n)<\infty$.  Fix $o\in M$, and set $X_n:= \{ (x_1,x_2)\in 
	 M\times M : n\leq \d(x_i,o)\leq n+1,\ \d(x_1,x_2)\leq r_n \}$, 
	 $Y_n:= \cup_{k=1}^n X_k$, $n\leq\infty$, and finally let $T_n$ be 
	 the integral operator whose kernel, a section of $End(F)$ denoted 
	 $k_n$, is the characteristic function of $Y_n$.  Since $k_n$ has 
	 compact support, if $n<\infty$, $\f(T_n)=0$.  On the contrary, 
	 since $Y_\infty$ contains the diagonal of $M\times M$, clearly 
	 $\f(T_\infty)=1$.  Finally
	 \begin{align*}
		 \|T_\infty-T_n\| & \leq \sup_{x\in M} \int_M 
		 \chi_{\cup_{k=n+1}^\infty X_k} (x,y)dy \\
		 & \leq \sup_{x\in M} \sum_{k=n+1}^\infty \int_M 
		 \chi_{X_k}(x,y)dy \\
		 & \leq \sup_{x\in M} \sum_{k=n+1}^\infty V(x,r_k) \\
		 & \leq \sum_{k=n+1}^\infty \b_2(r_k) \to0. 
	 \end{align*}
	 This proves both the assertions.
 \end{proof}
 
 \begin{Dfn}
	 Denote by $Tr_\ck$ the lower-semicontinuous semifinite trace on 
	 $\ca(F)$ obtained by regularising $\f|_\ca$, as in the previous 
	 subsection.
 \end{Dfn}

 \begin{Prop}\label{Prop:compzero}
	 $Tr_\ck$ vanishes on compact operators.
 \end{Prop}
 \begin{proof}
	 If $e$ is the one-dimensional projection onto the multiples of 
	 $f\in L^2(F)$, $\f(e) = \lo \frac{Tr(E_{K_{n}}eE_{K_{n}}}{vol 
	 K_{n}} = \lo \frac{\int_{K_{n}}|\f(x)|^{2}dx}{vol K_{n}} = 0$.  
	 Therefore $0\leq Tr_\ck(e) \leq \f(e) =0$, and $Tr_\ck(T)=0$ for 
	 any posite finite rank operator.  Let now $T$ be a compact 
	 operator, so that $T=U|T|$, and $|T|$ is the norm limit of a 
	 sequence $S_n$ of positive finite rank operators.  Then $0 \leq 
	 |Tr_\ck(T)| \leq \|U\| Tr_\ck(|T|) \leq \liminf Tr_\ck(S_n) =0$.
 \end{proof}

 Finally we give a sufficient criterion for a positive operator $A\in 
\ca$ to 
 satisfy $Tr_\ck(A)=\f(A)$.

 \begin{Prop}\label{3.4.2} 
	 Let $A\in\cj_{0+}$ be an integral operator, whose kernel $a(x,y)$ 
	 is a section of $End(F)$ which is uniformly continuous in a 
	 neighborhood of the diagonal in $M\times M$, namely
	 \begin{equation}\label{e:uniform}
		 \forall \eps>0,\ \exists \d>0\ : \d(x,y)<\d \imply 
		 |a(x,y)-a(x,x)|<\eps.
	 \end{equation}
	 Then $Tr_\ck(A)=\f(A)$.
 \end{Prop}
 \begin{proof}
	 Consider first a family of integral operators $B_\d$, with 
	 kernels, which are sections of $End(F)$, given by
	 \begin{equation*}
		 b_\d(x,y):= \frac{\b_1(\d)}{\b_2(\d)}\ 
		 \frac{\chi_{\D_\d}(x,y)}{ V(x,\d)},
	 \end{equation*}
	 where $\D_\d := \{ (x,y)\in M\times M : \d(x,y)<\d \}$.  Then 
	 $\sup_{x\in M} \int_M b_\d(x,y)dy = \frac{\b_1(\d)}{ \b_2(\d)} 
	 \leq 1$, and $\sup_{y\in M} \int_M b_\d(x,y)dx \leq 
	 \frac{\sup_{y\in M} V(y,\d)}{\b_2(\d)} \leq 1$, which imply 
	 $\|B_\d\|\leq1$, by Riesz-Thorin theorem.  \\
	 Set $E_n$ for the multiplication operator by the characteristic 
	 function of $K_{n}$, and observe that
	 \begin{align*}
		 Tr(E_nB_\d B_\d^*E_n) & = \int_{K_{n}}dx\int_M b_\d(x,y)^2dy 
		 \\
		 & = \frac{\b_{1}(\d)^{2}}{\b_{2}(\d)^{2}}\int_{K_{n}} 
		 \frac{dx}{V(x,\d)} \\
		 & \leq \frac{\b_1(\d)}{\b_2(\d)^2}\ vol(K_{n}) \leq 
		 \frac{vol(K_{n})}{\b_2(\d)}
	 \end{align*}
	 Therefore $\f(B_\d B_\d^*)\leq \b_2(\d)^{-1}$.  This implies 
	 that $\psi_\d := \f(B_\d\cdot B_\d^*)$ belongs to $\ca^*_+$, 
	 and $\psi_\d\leq \f|_\ca$ by Lemma \ref{3.2.8}.  By the results of 
	 the previous subsection, we have $\psi_\d(A) \leq Tr_\ck(A) \leq 
	 \f(A)$, for any $A\in\ca_+$.  \\
	 Take now $A\in\ca_+$ satisfying (\ref{e:uniform}), for a pair 
	 $\eps>0$, $\d>0$, and, setting $\b(\d):= (\frac{\b_1(\d) }{ 
	 \b_2(\d)})^2$ to improve readability, compute
	 \begin{align*}
		 |Tr & (E_nB_\d AB_\d^*E_n) - Tr(E_nAE_n)| \\
		 & \leq |Tr(E_nB_\d AB_\d^*E_n) - \b(\d) Tr(E_nAE_n)| + 
		 (1-\b(\d)) Tr(E_nAE_n) \\
		 & \leq \int_{K_{n}} dx \int_{B(x,\d)\times B(x,\d)} 
		 b_\d(x,y)|a(y,z)-a(x,x)|b_\d(x,z) dydz \\
		 & \qquad + (1-\b(\d)) Tr(E_nAE_n) \\
		 & \leq 3\eps \int_{K_{n}} dx \int_{B(x,\d)\times B(x,\d)} 
		 b_\d(x,y)b_\d(x,z) dydz \\
		 & \qquad + (1-\b(\d)) Tr(E_nAE_n) \\
		 & \leq 3\eps \b(\d) vol(K_{n}) + (1-\b(\d)) Tr(E_nAE_n) 
	 \end{align*}
	 By the arbitrariness of $\eps$ we get
     \begin{equation*}
		 \frac{|Tr(E_nB_\d AB_\d^*E_n) - Tr(E_nAE_n)|}{vol(K_{n})} 
		 \leq (1-\b(\d)) \frac{Tr(E_nAE_n)}{vol(K_{n})}.
     \end{equation*}
	 This implies $|\psi_\d(A)-\f(A)| \leq (1-\b(\d)) \f(A)$.  By 
	 Lemma \ref{2.1.2}$(ii)$, we get the thesis.
 \end{proof}

 \begin{Prop}\label{3.4.3} 
	 For any $t>0$, $\e{-t\D_p}$ belongs to the domain of $Tr_\ck$ and 
	 $Tr_\ck(\e{-t\D_{p}}) = \f(\e{-t\D_{p}})$, where $\D_{p}$ is the 
	 $p$-Laplacian operator.
 \end{Prop}
 \begin{proof}
	 By Remark~\ref{rem:laplacian} and Corollary~\ref{3.1.6} we have 
	 that $\e{-t\D_p}$ belongs to $\cj_{0+}\cap\ca$, hence, by 
	 Corollary~\ref{3.3.2}, $(iii)$, it belongs to the domain of 
$Tr_\ck$.  
	 The equality then follows by Proposition~\ref{3.4.2}.
 \end{proof}
 
 \begin{Prop}\label{prop:t0=t}
 	$Tr_\ck(f(\D_{p})) = \f(f(\D_{p}))$, for any $f\in C_{c}[0,\infty)$.
 \end{Prop}
 \begin{proof} 
	Let us introduce the positive functionals $T_{0}: f\in 
	C_{c}[0,\infty)_{+} \to \f(f(\D_{p}))\in [0,\infty)$, and $T: f\in 
	C_{c}[0,\infty)_{+} \to Tr_\ck(f(\D_{p}))\in [0,\infty)$.  Then, 
	by Riesz theorem, there are regular Borel measures $\m_{0}$, $\m$ 
	on $[0,\infty)$ s.t. $T_{0}(f) = \int f d\m_{0}$ and $T(f) = \int 
	f d\m$.  Then (see $e.g.$ \cite{GI4}, Propositions 5.2, 5.4) 
	$T_{0}(f) \geq \int f d\m_{0}$, $T(f) = \int f d\m$, for $f\in 
	C_{0}[0,\infty)_{+}$.  Let us prove that $\int f d\m_{0} \geq \int 
	f d\m$, for $f\in C_{0}[0,\infty)_{+}$.  Indeed, setting 
	$\widetilde{T_{0}}(f) := \int f d\m_{0}$, for $f\in 
	C_{0}[0,\infty)_{+}$, we get a semicontinuous weight on 
	$C_{0}[0,\infty)$, so that $\widetilde{T_{0}}(f) \equiv \sup \{ 
	S(f) : S\in C_{0}[0,\infty)^{*}_{+}, S\leq T_{0} \}$.  By 
	definition, $Tr_\ck(a) \equiv \sup \{ \psi(a) : \psi\in 
	\ca^{*}_{+}, \psi\leq \f|_\ca \}$, for any $a\in\ca_{+}$ for which 
	the right-hand side is finite.  Therefore for any $\psi\in 
	\ca^{*}_{+}$ s.t. $\psi\leq \f|_\ca$, setting $S(f) := 
	\psi(f(\D_{p}))$, we get $S\in C_{0}[0,\infty)^{*}_{+}$ and $S\leq 
	T_{0}$, so that $\psi(f(\D_{p})) \equiv S(f) \leq 
	\widetilde{T_{0}}(f)$, and, from the arbitrariness of $\psi$, 
	$T(f) = Tr_\ck(f(\D_{p})) \leq \widetilde{T_{0}}(f)$.  \\
	Therefore we conclude that $\m_{0}-\m$ is a positive measure on 
	$[0,\infty)$.  As $\f(\e{-t\D_{p}})=Tr_\ck(\e{-t\D_{p}})$, $t>0$, 
	we have $\int \e{-t\l} d(\m_{0}-\m)(\l) = 0$, so that $\m_{0}=\m$, 
	which implies $Tr_\ck(f(\D_{p})) = \f(f(\D_{p}))$, for $f\in 
	C_{c}[0,\infty)_{+}$.
 \end{proof}
 
 \section{Novikov-Shubin invariants and singular traces}
 \label{sec:singular}
 
 In this section we consider an open manifold with C$^{\infty}$ 
 bounded geometry possessing a regular exhaustion, i.e. the same 
 hypotheses assumed in subsection \ref{subsec:roe}.  Let us fix a 
 $p\in\{0,\ldots,n\}$ and denote by $\ca_p\equiv \ca(\La^{p}T^{*}M)$ 
 the C$^{*}$-algebra of almost local operators acting on 
 $L^{2}(\La^{p}T^{*}M)$, the Hilbert space of $L^{2}$-sections of the 
 vector bundle $\La^{p}T^{*}M$, $Tr_\ck$ denotes the 
 lower-semicontinuous semifinite trace on $\ca$ obtained in subsection 
 \ref{subsec:semicontinuous}, and $\D_{p}$ is the $p$-Laplacian, 
 acting as a selfadjoint operator on $L^{2}(\La^{p}T^{*}M)$.

 \subsection{Novikov-Shubin numbers for open manifolds and their 
 invariance} \label{subsec:laplacian}

 In this subsection we define the Novikov-Shubin numbers for such 
 manifolds and prove their invariance under quasi-isometries.  To do 
 this we  enlarge the C$^{*}$-algebra $\ca_p$ in order to 
 include sufficiently many spectral projections of the $p$-Laplacian.  
 It turns out that the noncommutative analogue of Riemann integrable 
 functions developed in \cite{GI4} gives a convenient framework, so we 
 first recall its construction and basic properties.
 
 Let $\ca$ be a general C$^{*}$-algebra with a semicontinuous 
 semifinite trace $\t$.  A pair of families $(A^-,A^+)$ in $\ca_{\sa}$ 
 is called an $\car$-cut in $\ca$ w.r.t. $\t$ if they are bounded, 
 separated (i.e. $a^-\leq a^+$ for any $a^\pm\in A^\pm$) and 
 $\t$-contiguous (i.e. $\forall\eps>0$ $\exists a^{\pm}_\eps\in A^\pm$ 
 s.t. $\t(a^+-a^-)<\eps)$.  A selfadjoint element $x\in\ca''$ is said 
 separating for $(A^-,A^+)$ if $a^-\leq x\leq a^+$ for any $a^\pm\in 
 A^\pm$.  Then the following holds.
 
 \begin{Thm} {\rm \cite{GI4}} 
	 The set of separating elements between $\car$-cuts in $\ca$ is 
	 the selfadjoint part of a C$^*$-algebra, denoted by $\ar$, and 
	 called the C$^*$-algebra of Riemann measurable operators.  The 
	 GNS representation $\pi_\tau$ of $\ca$ extends to a 
	 $^*$-homomorphism (still denoted by $\pi_\tau$) of $\ar$ to 
	 $\pi_\tau(\ca)''$, hence the trace $\t$ extends to a 
	 semicontinuous semifinite trace on $\ar$.  $\ar$ contains all the 
	 separating elements between $\car$-cuts in it, and is closed 
	 under functional calculus with Riemann integrable functions, in 
	 particular almost all spectral projections $e_{[a,b)}$ of 
	 selfadjoint elements of $\ca$ belong to $\ar$.
 \end{Thm}

 Applying this result to $\ca_p$, we obtain the C$^*$-algebra 
 $\ca^\car_p$ with a lower-semicontinuous semifinite trace, still 
 denoted $Tr_\ck$.  Then $\chi_{[0,t)}(\D_p)$ and 
 $\chi_{[\eps,t)}(\D_p)$ are Riemann measurable spectral projections 
 and belong to $\ar_p$ for almost all $t>\eps>0$, by the previous 
 Theorem.  Denote by $N_{p}(t) := Tr_\ck(\chi_{[0,t)}(\D_p))$, 
 $\th_{p}(t):=Tr_\ck(\e{-t\D_p})$.

 \begin{Lemma}
	$\th_p(t) = \int_{0}^{\infty}\e{-t\l} dN_p(\l)$ so that 
	$\lim_{t\to0} N_p(t)= \lim_{t\to\infty}\th_p(t)$.
 \end{Lemma}
 \begin{proof}
	If $\D=\int_{0}^{\infty}\l de(\l)$ denotes the spectral 
	decomposition, then $\e{-t\D}=\int_{0}^{\infty}\e{-t\l} de(\l)$.  
	Since the latter is defined as the norm limit of the 
	Riemann-Stieltjes sums, $\pi_p(\e{-t\D}) = 
	\int_{0}^{\infty}\e{-t\l} d\pi_p(e(\l))$, where $\pi_p$ denotes 
	the GNS representation of $\ca_p$ w.r.t the trace $Tr_\ck$.  The 
	result then follows by the normality of the trace in the GNS 
	representation.
 \end{proof}
 
 \begin{Dfn} \label{def:NS-inv} 
	 We define $b_{p}\equiv b_p(M,\ck) :=\lim_{t\to0} N_{p}(t)= 
	 \lim_{t\to\infty}\th_{p}(t)$ to be the $p$-th L$^{2}$-Betti 
	 number of the open manifold $M$ endowed with the exhaustion 
	 $\ck$.  Let us now set $N^{0}_{p}(t) := N_{p}(t) - b_{p} \equiv 
	 \lim_{\eps\to0} Tr_\ck(\chi_{[\eps,t)}(\D_p))$, and 
	 $\th^{0}_{p}(t) := \th_{p}(t) - b_{p} = \int_{0}^{\infty}\e{-t\l} 
	 dN^{0}_{p}(\l)$.  The Novikov-Shubin numbers of $(M,\ck)$ are 
	 then defined as
	 \begin{align*}
		 \a_{p}\equiv \a_{p}(M,\ck) & := 2\limsup_{t\to0} \frac{\log 
		 N^{0}_{p}(t)}{\log t}, \\
		 \underline{\a}_{p}\equiv \underline{\a}_{p}(M,\ck) & := 
		 2\liminf_{t\to0} \frac{\log N^{0}_{p}(t)}{\log t}, \\
		 \a'_{p}\equiv \a'_{p}(M,\ck) & := 2\limsup_{t\to\infty} 
		 \frac{\log \th^{0}_{p}(t)}{\log 1/t}, \\
		 \underline{\a}'_{p}\equiv \underline{\a}'_{p}(M,\ck) & := 
		 2\liminf_{t\to\infty} \frac{\log \th^{0}_{p}(t)}{\log 1/t}.
	 \end{align*}
 \end{Dfn}
 
 It follows from (\cite{GS}, Appendix) that $\underline{\a}_{p} = 
 \underline{\a}'_{p}\leq \a'_{p}\leq \a_{p}$, and 
 $\a'_{p} = \a_{p}$ if $\th_{p}^{0}(t)= O(t^{-\d})$, for 
 $t\to\infty$, or equivalently $N^{0}_{p}(t)=O(t^{\d})$, for $t\to0$.
 Observe that L$^{2}$-Betti numbers and Novikov-Shubin numbers depend 
 on the limit procedure $\om$ and the exhaustion $\ck$.

 \begin{Prop}\label{p:a0=a'0}
	Let $M$ be a complete non-compact Riemannian manifold of positive 
	injectivity radius and Ricci curvature bounded from below.  Then 
	$\a_{0}(M,\ck)=\a'_{0}(M,\ck)\geq 1$ for any regular exhaustion 
	$\ck$.
 \end{Prop}
 \begin{proof}
	 Recall that, under the previous assumptions, Varopoulos 
	 \cite{Varopoulos2} proved that the heat kernel on the diagonal 
	 has a uniform inverse-polynomial bound, more precisely, in the 
	 strongest form due to \cite{ChavelFeldman}, we have
	 $$
	 \sup_{x,y\in M} H_{0}(t,x,y)\leq Ct^{-1/2}
	 $$
	 for a suitable constant $C$.  Then, as 
	 $$
	 \th(t)=\t(\e{-t\D})= \lo 
	 \frac{\int_{B(o,n_{k})}H_{0}(t,x,x)dvol(x)}{V(o,n_{k})} \leq 
	 Ct^{-1/2},
	 $$
	 it follows from (\cite{GS}, Appendix) that 
	 $\a_{0}=\a_{0}'$, which concludes the proof.
 \end{proof}
  
 \begin{rem}\label{r:remarks}
	$(a)$ If $M$ is a covering of a compact manifold $X$, L$^{2}$-Betti 
	numbers were introduced by Atiyah \cite{Atiyah} whereas 
	Novikov-Shubin numbers were introduced in \cite{NS1}.  They were 
	proved to be $\G$-homotopy invariants, where $\G:=\pi_{1}(X)$ is 
	the fundamental group of $X$, by Dodziuk \cite{Dodziuk} and 
	Gromov-Shubin \cite{GS} respectively.  L$^{2}$-Betti numbers were 
	subsequently defined for the open manifolds considered in this 
	paper by Roe \cite{Roe1}, though in a different way, and were 
	proved to be invariant under quasi-isometries (see below) in 
	\cite{RoeBetti}.
	\\
	$(b)$ In the case of coverings, the trace $Tr_\Gamma$ is
	normal on the von~Neumann algebra of $\G$-invariant operators,
	hence $\lim_{t\to0}Tr(e_{[0,t)}(\D_p))=Tr(e_{\{0\}}(\D_p))$.  In
	the case of open manifolds there is no natural von~Neumann
	algebra containing the bounded functional calculi of $\D_p$ on
	which the trace $Tr_\ck$ is normal, hence the previous
	equality does not necessarily hold. Such phenomenon has been
	considered by Farber in \cite{Farber} in a context which is
	similar to ours, and the difference
	$\lim_{t\to0}Tr(e_{[0,t)}(\D_p))-Tr(e_{\{0\}}(\D_p))$ has been
	called the {\it torsion dimension}. We shall denote by
	$\tordim(M,\D_{p})$ such difference, and shall sometimes assume it
	vanishes. We are not aware of a general vanishing result in
	our context.
	\\
	$(c)$ Let us observe that the above definitions for
	$L^2$-Betti numbers and Novikov-Shubin numbers coincide with
	the classical ones in the case of amenable coverings, if one
	chooses the exhaustion given by the F\o lner condition.  An
	explicit argument is given in \cite{GI5}.
	\\
	$(d)$ In the case of coverings there is a well-known conjecture on 
	the positivity of the $\a_{p}$'s.  A result by Varopoulos 
	\cite{Varopoulos1} shows that $\a_{0}$ is a positive integer, hence 
	$\a_{0}\geq1$.  The previous Proposition extends this inequality 
	to the case of open manifolds.  Moreover $\a_{0}$ can assume any 
	value in $[1,\infty)$, as follows from Theorem \ref{a0=dinfinity} 
	and Corollary \ref{2.3.4}.
 \end{rem}

 Our first objective is to show that our definition of L$^{2}$-Betti 
 numbers coincides with Roe's definition.  Then we prove that 
 Novikov-Shubin numbers are invariant under quasi-isometries, where a 
 map $\f: M\to \widetilde M$ between open manifolds of \bg is a 
 quasi-isometry \cite{RoeBetti} if $\f$ is a diffeomorphism s.t.
 \begin{itemize}
	 \itm{i} there are $C_{1},\ C_{2}>0$ s.t. $C_{1}\|v\| \leq 
	 \|\f_{*}v\| \leq C_{2}\|v\|$, $v\in TM$ 
	 \itm{ii} $\nabla-\f^{*}\widetilde\nabla$ is bounded with all its 
	 covariant derivatives, where $\nabla$, $\widetilde\nabla$ are 
	 Levi-Civita connections of $M$ and $\widetilde M$.
 \end{itemize}
 
 Recall Roe's definition of L$^{2}$-Betti numbers \cite{Roe2}, which 
 we temporarily denote by $b_{p}^{R}$.  $b_{p}^{R} := \inf \{ 
 \phi(f(\D_{p})) ; f\in\cc \} = \inf \{Tr_\ck(f(\D_{p})) ; f\in\cc 
 \}$, by Proposition \ref{prop:t0=t}, where $\cc := \{ f \in 
 C_{c}^{\infty}[0,\infty)_{+} : f(0)=1 \}$.  Then
 
 \begin{Prop}
 	$b^{R}_{p} = b_{p}$.
 \end{Prop}
 \begin{proof}
	As $b_{p} \equiv \inf_{t>0} Tr_\ck(\chi_{[0,t)}(\D_{p}))$, and for 
	any $t>0$ there is $f\in\cc$ s.t. $f\leq \chi_{[0,t)}$, we 
	conclude that $b_{p}^{R} \leq b_{p}$.  Let now $\eps>0$ be given, 
	and let $f\in\cc$ be s.t. $(1+\eps) b_{p}^{R} > 
	Tr_\ck(f(\D_{p}))$.  Then there is $\d=\d_{\eps}>0$ s.t. $f(x) 
	\geq 1-\eps$, for $x\in[0,\d)$, so that $\frac{1}{1-\eps} f \geq 
	\chi_{[0,\d)}$ and $(1+\eps) b_{p}^{R} \geq (1-\eps) 
	Tr_\ck(\chi_{[0,\d)}(\D_{p})) \geq (1-\eps) b_{p}$, that is 
	$b^{R}_{p} \geq \frac{1-\eps}{1+\eps} b_{p}$, and from the 
	arbitrariness of $\eps$, $b^{R}_{p}\geq b_{p}$, and the thesis 
	follows.
 \end{proof}
  
 \begin{Thm} {\rm \cite{RoeBetti}} 
	 $L^2$-Betti numbers are invariant under quasi-isometries.
 \end{Thm}

 \begin{Thm} \label{thm:invariance} 
	 Let $(M,\ck)$ be an open manifold with a regular exhaustion, and 
	 let $\f:M\to\widetilde M$ be a quasi-isometry.  Then $\f(\ck)$ is 
	 a regular exhaustion for $\widetilde M$, 
	 $\a_{p}(M,\ck)=\a_{p}(\widetilde M,\f(\ck))$ and the same holds 
	 for $\underline{\a}_{p}$ and $\a'_{p}$.
 \end{Thm}
 \begin{proof} 
	  Let us denote by 
	  $\Phi\in\cb(L^{2}(\La^{*}T^{*}M),L^{2}(\La^{*}T^{*}\widetilde 
	  M))$ the extension of $(\f^{-1})^{*}$.  Then 
	  $Tr_{\f(\ck)}=Tr_\ck(\Phi^{-1}\cdot\Phi)$.  Also, setting 
	  $e_{\eps,t}:=\chi_{[\eps,t)}(\D_p)$, 
	  $q_{\eta,s}:=\Phi^{-1}\chi_{[\eta,s)}(\widetilde\D_p)\Phi$, we 
	  have
	 \begin{align*}
		0 & \leq Tr_\ck(e_{\eps,t}-e_{\eps,t}q_{\eta,s}e_{\eps,t}) = 
		Tr_\ck(e_{\eps,t}(1-q_{\eta,s})e_{\eps,t}) \\
		& = Tr_\ck(e_{\eps,t}q_{0,\eta}e_{\eps,t}) + 
		Tr_\ck(e_{\eps,t}q_{s,\infty}e_{\eps,t}) \\
		& = Tr_\ck(q_{0,\eta}e_{\eps,t}q_{0,\eta}) + 
		Tr_\ck(e_{\eps,t}e_{0,t}q_{s,\infty}e_{0,t}) \\
		& \leq Tr_\ck(q_{0,\eta}e_{\eps,\infty}q_{0,\eta}) + 
		Tr_\ck(e_{\eps,t}e_{0,t}q_{s,\infty}e_{0,t}) \\
		& \leq Tr_\ck(q_{0,\eta})\|q_{0,\eta}e_{\eps,\infty}
		q_{0,\eta}\| + 
		Tr_\ck(e_{\eps,t})\|e_{0,t}q_{s,\infty}e_{0,t}\| \\
		& \leq Tr_\ck(q_{0,\eta})\ C \sqrt{\frac{\eta}{\eps}} + 
		Tr_\ck(e_{\eps,t})\ C \sqrt{\frac t s},
	\end{align*}
	where the last inequality follows from \cite{RoeBetti}.
	Then
	\begin{align*}
		Tr_\ck(q_{\eta,s}) & = Tr_\ck(e_{\eps,t}) + 
		Tr_\ck(q_{\eta,s}-e_{\eps,t}q_{\eta,s}e_{\eps,t})- 
		Tr_\ck(e_{\eps,t}-e_{\eps,t}q_{\eta,s}e_{\eps,t}) \\
		& \geq Tr_\ck(e_{\eps,t}) - Tr_\ck(q_{0,\eta})\ C 
		\sqrt{\frac{\eta}{\eps}}-Tr_\ck(e_{\eps,t})\ C\sqrt{\frac t s}.
	\end{align*}
	Now let $a>1$ and compute
	\begin{align*}
		\widetilde N^{0}(s) 
		& = \lim_{\eps\to0}Tr_\ck(q_{\eps^{a},s}) \geq 
		\lim_{\eps\to0} \left[ Tr_\ck(e_{\eps,t}) - 
		Tr_\ck(q_{0,\eps^{a}})\ C 
		\eps^{\frac{a-1}2} - Tr_\ck(e_{\eps,t})\ C \sqrt{\frac t s} 
		\right] \\
		&  = N^{0}(t)\left[ 1-C\sqrt{\frac t s} \right].
	\end{align*}
	Therefore with $\l:=4C^{2}$ we get $\widetilde N^{0}(\l t) \geq 
	\frac12 N^{0}(t)$, and exchanging the roles of $M$ and $\widetilde 
	M$, we obtain $\frac12 N^{0}(\l^{-1}t)\leq \widetilde N^{0}(t) 
	\leq 2N^{0}(\l t)$.  This means that $N^{0}$ and $\widetilde 
	N^{0}$ are dilatation-equivalent (see \cite{GS}) so that the 
	thesis follows from \cite{GS}.
 \end{proof}

 \begin{rem}
	 We have chosen Lott's normalization \cite{Lott} for the 
	 Novikov-Shubin numbers $\a_{p}(M)$ because Laplace operator is a 
	 second order differential operator, and this normalization gives 
	 the equality between $\a_{0}(M)$ and the asymptotic dimension of 
	 $M$, cf.  Theorem \ref{a0=dinfinity}.  \\ Our choice of the 
	 $\limsup$ in Definition \ref{def:NS-inv}, in contrast with 
	 Lott's choice \cite{Lott}, is motivated by our interpretation 
	 of $\a_{p}(M)$ as a dimension.  On the one hand, a noncommutative 
	 measure corresponds to $\a_p$ via a singular trace, according to 
	 Theorem \ref{thm:ap}.  On the other hand, this choice implies 
	 that $\a_0$, being equal to the asymptotic dimension of $M$, 
	 possesses the classical properties of a dimension as stated in 
	 Theorem \ref{Thm:adim}, cf.  also Remark \ref{1.2.7} $(b)$.
 \end{rem}

\subsection{Novikov-Shubin numbers as asymptotic spectral dimensions} 
\label{subsec:singtrac}
 In this subsection we discuss a dimensional interpretation for the 
 Novikov-Shubin numbers.  In the case of compact manifolds, the 
 dimension may be recovered by the Weyl asymptotics.  In particular, 
 the formula $(\lim\frac{\log\m_n}{\log 1/n})^{-1}$ gives the 
 dimension of the manifold, where $\m_{n}$ refers to $\D^{-1/2}$.  
 This formula makes sense also in the non-compact case, if one 
 replaces the eigenvalue sequence with the eigenvalue function as 
 explained below, and clearly recovers the dimension of the manifold.  
 But in this case, the behaviour for $t\to0$ may be considered also, 
 giving rise to an asymptotic counterpart of the dimension.  In the 
 case of covering manifolds \cite{GI5}, we defined the {\it asymptotic 
 spectral dimension of the triple $(M,\G,\D_p)$} as 
 $\left(\liminf_{t\to0}\frac{\log\m_{p}(t)}{\log1/t}\right)^{-1}$, 
 where $\m_{p}$ refers to the operator $\D_p^{-1/2}$.

 On the one hand, this number is easily shown to coincide with the
 $p$-th Novikov-Shubin number; on the other hand, it deserves the name
 of dimension also in the context of noncommutative measure
 theory. Indeed, Hausdorff dimension determines which power of the
 radius of a ball gives rise to a non trivial volume on the
 space. Analogously, the spectral dimension determines which power of
 the $p$-Laplacian gives rise to a non trivial singular trace on the
 algebra $\ca_{p}$.

 We extend this result to the case of open manifolds, using the
 unbounded Riemann integration and the theory of singular traces for
 C$^*$-algebras developed in \cite{GI4}. However, since the trace we
 use is not normal with respect to the given representation of $\ca_{p}$
 on the space of $L^2$-differential forms, some assumptions like the
 vanishing of the torsion dimension introduced in Remark
 \ref{r:remarks} $(b)$ are needed.

 \medskip

 Let us briefly recall the definition and main properties of unbounded
 Riemann integrable operators.  A linear operator $T$ on
 $\ch$ is said to be {\it affiliated} to a von~Neumann algebra $\cam$
 ($T\aff \cam$) if all elements of $x\in\cam'$ send its domain into
 itself and $Tx\eta=xT\eta$, for any $\eta$ in $\cd(T)$.

 Let $\ca\subset\cb(\ch)$ be a C$^{*}$-algebra with a lower 
 semicontinuous semifinite trace $\t$, and $\ar$ be the 
 C$^{*}$-algebra of Riemann measurable operators.  A sequence 
 $\{e_n\}$ of projections in $\cb(\ch)$ is called a {\it Strongly 
 Dense Domain} (SDD) w.r.t. $\ar$ if $e_n^{\perp}\in\ar$ is 
 $\tau$-finite and $\tau(e_n^{\perp})\to0$.  We shall denote by $e$ 
 the projection $\sup_ne_n$.  Let us remark that, if the trace $\tau$ 
 is not faithful, $e$ is not necessarily $1$.  Nevertheless it is easy 
 to show that $e^{\perp}\in\ar$ and $\tau(e)=0$.

 \begin{Dfn} 
	 We denote by $\ov\ar$ the family of closed, densely defined 
	 operators affiliated to $\ca''$ for which there exists a SDD 
	 $\{e_n\}$ such that 
	 \itm{i} $e_n\ch\subset\cd(T)\cap\cd(T^{*})$  
	 \itm{ii} $eTe_n\in\ar$, $e_nTe\in\ar$.
 \end{Dfn}

 We also introduce the relation of $\t$-a.e. equality, which 
 turns out to be an equivalence relation, among operators in 
 $\ov\ar$, namely $T$ is equal to $S$ $\t$-a.e. if there 
 exists a common SDD $\{e_n\}$ for $T$ and $S$ such that 
 $eTe_{n}=eSe_{n}$ for any $n\in\bn$.  

 In the following we shall denote by $\pi$ the GNS representation of
 $\ca$ associated with the trace $\tau$, by $\cam$ the von~Neumann
 algebra $\pi(\ca)''$, and by $\mt$ the algebra of $\t$-measurable
 operators affiliated to $\cam$.

 \begin{Thm} {\rm \cite{GI4}} The set $\ov\ar$ is a
	 $\t$-a.e. bimodule, namely it is closed under strong sense
	 operations, and the usual properties of a $^*$-bimodule over
	 $\ar$ hold $\tau$-almost everywhere.  Moreover the GNS
	 representation extends to a map from $\ov\ar$ to $\mt$ which
	 preserves strong sense operations, hence the trace $\t$
	 extends to $\ov\ar$.  $\ov\ar$ contains the functional
	 calculi of the selfadjoint elements in $\ar$ under unbounded
	 Riemann integrable functions, and is called the $\t$-a.e.
	 bimodule of unbounded Riemann integrable elements. 
 \end{Thm}

 We define the distribution function w.r.t. $\t$ of an operator 
 $A\in\ov\ar$ via the distribution function of its image under the GNS 
 representation, namely $\l_A=\l_{\pi(A)}$.  If $A\in\ov\ar$ is a 
 positive (unbounded) continuous functional calculus of an element in 
 $\ca$, then $\c_{(t,+\infty)}(A)$ belongs to $\ar$ a.e., therefore 
 its distribution function is given by 
 $\l_A(t)=\t(\c_{(t,+\infty)}(A))$.  The non increasing rearrangement 
 is defined as $\m_A(t):=\inf\{s\geq0:\l_A(s)\leq t\}$.

 Now we come back to our concrete situation, namely to the pair
 $(\ca_p,Tr_\ck)$ associated with a manifold $M$ endowed with a regular
 exhaustion $\ck$. In the following, when the Laplacian $\D_p$ has a
 non trivial kernel, we denote by $\D_p^{-\a}$, $\a>0$, the
 (unbounded) functional calculus of $\D_p$ w.r.t. the function $\f_\a$
 given by $\f_\a(0)=0$ and $\f_\a(t)=t^{-\a}$ when $t>0$.

 \begin{Lemma}\label{Lemma:tordimass}
	 \itm{a} The following are equivalent 
	 \begin{itemize}
		 \itm{a.i} The projection $E_{p}$ onto the kernel of $\D_{p}$
		 (also expressed by $\chi_{\{1\}}(\e{-\D_{p}})$) is
		 Riemann measurable, and the torsion dimension
		 vanishes, namely $Tr_\ck(E_{p})$
		 is equal to $b_p$.
		 \itm{a.ii} $\chi_{\{1\}}(\pi(\e{-\D_{p}}))$ is
		 Riemann integrable in the GNS representation $\pi$.
	 \end{itemize}
	 \itm{b} The vanishing of the Betti number $b_p$ implies
	 $(a)$. It is equivalent to $(a)$ if $\ker(\D_{p})$ is
	 finite-dimensional.
	 \itm{c} If $(a)$ is satisfied, $\D_{p}^{-\a}\in\ov{\ar}$ for
	 any $\a>0$.
 \end{Lemma}
 \begin{rem}
 	The extension of the GNS representation $\pi$ to $\ar_{p}$ does not 
 	necessarily commute with the Borel functional calculus. In 
 	particular $\chi_{\{1\}}(\pi(\e{-\D_{p}}))$ is not necessarily equal 
 	to $\pi(\chi_{\{1\}}(\e{-\D_{p}}))$.
 \end{rem}
 \begin{proof}
	 $(a)$. $(i)\imply (ii)$. Since the projection
	 $E_{p}\equiv\chi_{\{0\}}(\D_p)$ is Riemann integrable and less than
	 $\e{-t\D_p}$ for any $t$, its image in the GNS
	 representation is Riemann integrable and less than
	 $\pi(\e{-t\D_p})$ for any $t$. This implies that
	 $\pi(E_p)\leq\chi_{\{1\}}(\pi(\e{-\D_{p}}))\leq \pi(\e{-t\D_p})$
	 is an $\car$-cut, hence the thesis.
	 \\
	 $(ii) \imply (i)$.  By normality of the trace in the GNS
	 representation, $Tr_\ck(\e{-t\D_p})$ converges to
	 $Tr_\ck(\chi_{\{1\}}(\pi(\e{-\D_{p}})))$ hence, by hypothesis,
	 for any $\eps>0$ we may find $a_\eps\in\ca$ and $t_\eps>0$
	 s.t. $a_\eps\leq\chi_{\{1\}}(\pi(\e{-\D_{p}}))$ and
	 $Tr_\ck(\e{-t_\eps\D_{p}} - a_\eps)<\eps$.  This implies
	 $a_\eps\leq\e{-t\D_p}$ for any $t$, hence
	 $a_\eps\leq\chi_{\{0\}}(\D_p)$, which means that
	 $(\{a_\eps\},\{\e{-t_\eps\D_p}\})$ is an $\car$-cut for
	 $\chi_{\{0\}}(\D_p)$, namely this projection is Riemann
	 integrable and
	 $Tr_\ck(\e{-t_\eps\D_p}-\chi_{\{0\}}(\D_p))\leq\eps$, i.e. the
	 thesis.
	 \\
	 $(b)$ If $\lim_{t\to\infty} Tr_\ck(\e{-t\D_{p}})=0$, from $0\leq
	 \chi_{\{1\}}(\e{-\D_{p}}) \leq \e{-t\D_{p}}$, we have that
	 $\chi_{\{1\}}(\e{-\D_{p}})$ is a separating element for an
	 $\car$-cut in $\ca$. The second statement follows from
	 Proposition~\ref{Prop:compzero}.
	 \\
	 $(c)$ We have to exhibit an SDD for $\D_{p}^{-\a}$.  Indeed,
	 as $\chi_{[t,\infty)}(\D_{p})$ is Riemann measurable for
	 almost all $t>0$, choose a strictly decreasing sequence
	 $t_{n}\to 0$ of such $t$, and set $e_{n}:=
	 \chi_{[t_{n},\infty)}(\D_{p})+ E_{p}$.
	 Then $e_{n}^{\perp}\equiv \chi_{(0,t_{n})}(\D_{p})\in\ar$,
	 $Tr_\ck(e_{n}^{\perp})= N^{0}_{p}(t_{n}) \to 0$, and
	 $e\D_{p}^{-\a}e_{n} \equiv \int_{t_{n}}^{\infty} \l^{-\a}
	 de(\l)\in \ar$.
 \end{proof}

 If hypothesis $(a)$ of the previous Lemma is satisfied, we may define 
 the distribution function $\l_{p}$ and the eigenvalue function 
 $\m_{p}$ for the operator $\D_p^{-1/2}$, hence the local spectral 
 dimension as the inverse of 
 $\lim_{t\to\infty}\frac{\log\m_{p}(t)}{\log1/t}$, and it may be shown 
 that such dimension gives the dimension of the manifold for any $p$.

 But, in the case of open manifolds, we may also consider the
 asymptotics for $t\to0$, which gives rise to an
 asymptotic counterpart of the dimension. Then we define the {\it
 asymptotic spectral dimension of the triple $(M,\ck,\D_p)$} as
 $$
 \left(\liminf_{t\to0}\frac{\log\m_{p}(t)}{\log1/t}\right)^{-1}.
 $$

 It is not difficult to show that such asymptotic spectral dimensions
 coincide with the Novikov-Shubin numbers.

 \begin{Thm}\label{Thm:singular} 
	 Let $(M,\ck)$ be an open manifold equipped with a regular 
	 exhaustion such that the projection on the kernel of $\D_p$ is 
	 Riemann integrable and  $\tordim(M,\D_{p})=0$.  
	 Then the asymptotic spectral dimension of $(M,\ck,\D_p)$ 
	 coincides with the Novikov-Shubin number $\a_p(M,\ck)$.
 \end{Thm}
 \begin{proof} 
	 By hypothesis, $e_{(0,t)}(\D_p)$ is Riemann integrable 
	 $Tr_\ck$-a.e., hence $N^0_p(t) = Tr_\ck(e_{(0,t)}(\D_p)) = 
	 Tr_\ck(e_{(t^{-1},\infty)}(\D_p^{-1})) = 
	 Tr_\ck(e_{(t^{-1/2},\infty)}(\D_p^{-1/2})) = \l_p(t^{-1/2})$.  Then
	 \begin{align}
		 \a_p &=2\limsup_{s\to0} \frac{\log N^0_{p}(s)}{\log{s}} 
		 &=2\limsup_{s\to0}\frac{\log\l_p(s^{-1/2})}{\log{s}} 
		 &=\limsup_{t\to\infty}\frac{\log\l_p(t)}{\log\frac1t}.
	 \end{align}
	 The statement follows if we prove that
	 $$
	 \liminf_{t\to0}\frac {\log\m(t)}{\log\frac1{t}} 
	 =\left(\limsup_{s\to\infty}\frac 
	 {\log\l(s)}{\log\frac1{s}}\right)^{-1}
	 $$
	 for any $t$, where the values 0 and $\infty$ are allowed.  
	 \\
	 Let $t_n\to0$ be a sequence such that 
	 $\lim_{n\to\infty}\frac{\log\m(t_n)}{\log\frac1{t_n}}= 
	 \liminf_{t\to0}\frac{\log\m(t)}{\log\frac1{t}}$, and let 
	 $t'_n:=\inf\{s\geq0:\m(s)=\m(t_n)\} 
	 =\min\{s\geq0:\m(s)=\m(t_n)\}$ where the last equality holds 
	 because of right continuity.  Then
	 $$
	 \liminf_{t\to0}\frac{\log\m(t)}{\log\frac1{t}} 
	 \leq\lim_{n\to\infty}\frac{\log\m(t'_n)}{\log\frac1{t'_n}} 
	 \leq\lim_{n\to\infty}\frac{\log\m(t_n)}{\log\frac1{t_n}} 
	 =\liminf_{t\to0}\frac{\log\m(t)} {\log\frac{1}{t}}
	 $$
	 namely we may replace $t_n$ with $t'_n$ to reach the $\liminf$.  
	 Also, $\l(\m(t'_n))=\inf\{t\geq0:\m(t)\leq\m(t'_n)\}=t'_n$, 
	 therefore
	 \begin{align*}
		 \liminf_{t\to0}\frac{\log\m(t)}{\log\frac{1}{ t}} 
		 &=\lim_{n\to\infty} \frac{\log\m(t'_n)}{\log\frac{1}{ t'_n}} 
		 =\lim_{n\to\infty} \frac{\log\m(t'_n)}{\log\frac{1}{ 
		 \l(\m(t'_n))}}\\
		 &\geq\liminf_{s\to0}\frac{\log s}{\log\frac{1}{ \l(s)}} 
		 =\left(\limsup_{s\to0}\frac{\log\l(s)}{\log\frac{1}{ 
		 s}}\right)^{-1}
	 \end{align*}
	 For the converse inequality, let $s_n\to\infty$ be a sequence for 
	 which $\lim_{n\to\infty} \frac{\log\l(s_n)}{\log\frac{1}{ 
	 s_n}}=\limsup_{s\to0} \frac{\log\l(s)}{\log\frac{1}{ s}}$.  As 
	 before, $s'_n:=\inf\{s\geq0:\l(s)=\l(s_n)\} 
	 =\min\{s\geq0:\l(s)=\l(s_n)\}$ still brings to the $\limsup$ and 
	 verifies $\m(\l(s'_n))=s'_n$.  
 \end{proof}

 Now we generalise a result proved in \cite{GI5} in the case of
 coverings, proving that the spectral dimensions defined above
 give rise to singular traces, namely select the
 correct power of the $p$-Laplace operator which gives rise to a non
 trivial singular trace on the $Tr_\ck$-a.e.  bimodule $\ov{\ca_{p}^\car}$.

 Let us recall that an operator $T\in\ov\ar$ is called $0$-{\it
 eccentric} when
 \begin{align*}
	\limsup_{t\to0}\frac{\int_0^t\m_T(s)ds}{\int_0^{2t}\m_T(s)ds}=1, & 
	\quad\mathrm{ if}\quad \int_0^1\mu_T(t)dt <\infty, \\
	\liminf_{t\to0}\frac{\int_t^1\m_T(s)ds}{\int_{2t}^1\m_T(s)ds}=1, & 
	\quad\mathrm{ if}\quad \int_0^1\mu_T(t)dt =\infty.
 \end{align*}
 
 As in the case of von~Neumann algebras, any $0$-eccentric operator in
 $\ov\ar$ gives rise to a singular trace, more precisely to a trace on
 $\ov\ar$ which vanishes on all bounded operators.  Singular traces
 may be described as the pull-back of the singular traces on $\cam$
 via the (extended) GNS representation.  On the other hand, explicit
 formulas may be written in terms of the non-increasing rearrangement.
 We write these formulas for the sake of completeness.

 \begin{Thm}\label{Thm:singtrac} {\rm \cite{GI4}}
	 If $T\in\ov\ar$ is $0$-eccentric and $\int_0^1\mu_T(t)dt 
	 <\infty$, there exists a generalized limit $\Lim_\om$ in $0$ such 
	 that the functional
	 $$ 
	 \t_\om(A):=\Lim_\om\left(\frac{\int_0^t\m_A(s)ds}{\int_0^{t}\m_T(s)ds}\right) 
	 \quad A\in X(T)_+
	 $$
	 extends to a singular trace on the a.e. $^*$-bimodule $X(T)$ over 
	 $\ar$ generated by $T$.  If $\int_0^1\mu_T(t)dt =\infty$, the 
	 previous formula should be replaced by
	 $$ 
	 \t_\om(A):=\Lim_\om\left(\frac{\int_t^1\m_A(s)ds}{\int_{t}^1\m_T(s)ds}\right), 
	 \quad A\in X(T)_+.
	 $$
	 Such traces naturally extend to traces on $X(T)+\ar$.
 \end{Thm} 

 Making use of the previous Theorem, we show that Novikov-Shubin
 numbers are dimensions in the sense of noncommutative measure theory.

 \begin{Thm}\label{thm:ap} 
	 Let $(M,\ck)$ be an open manifold with a regular exhaustion such 
	 that the projection on the kernel of $\D_p$ is Riemann integrable 
	 and $\tordim(M,\D_{p})=0$.  If $\a_p$ is finite 
	 nonzero, then $\D_p^{-\a_p/2}$ is $0$-eccentric, namely gives rise 
	 to a non trivial singular trace on $\ov{\ca^\car}$.
 \end{Thm}
 \begin{proof} 
	 The $0$-eccentricity of $\D_p^{-\a_p/2}$ follows by \cite{GI5}, 
	 hence the thesis follows by Theorem~\ref{Thm:singtrac}.
 \end{proof}

 \section{An asymptotic dimension for metric spaces.}\label{sec:first}

 The main purposes of this section are to introduce  an asymptotic 
 dimension for metric spaces, and to show that for a suitable 
 class of manifolds ($i.e.$ open manifolds of C$^{\infty}$-bounded 
 geometry satisfying Grigor'yan's isoperimetric inequality) the 
 asymptotic dimension coincides with the $0$-th Novikov-Shubin 
 invariant. In particular this shows that, for these manifolds, 
 $\a_{0}$ does not depend on the limit procedure $\om$ and, in a mild 
 sense, is independent of the regular exhaustion too. More 
 precisely, $\a_{0}$ does not change if we choose $\ck$ among the 
 exhaustions by balls with a common centre.

 To our knowledge, the notion of asymptotic dimension in the general 
 setting of metric dimension theory has not been studied, even though 
 Davies \cite{Davies} proposed a definition in the case of cylindrical 
 ends of a Riemannian manifold.  We shall give a definition of 
 asymptotic dimension for a general metric space, based on the (local) 
 Kolmogorov dimension \cite{KT} and state its main properties.  We 
 compare our definition with Davies'.

 \subsection{Kolmogorov-Tihomirov metric dimension} 
 \label{subsec:KTdim}

 In this subsection we recall a definition of metric dimension due to 
 Kolmogorov and Tihomirov \cite{KT} (see also \cite{Falconer} where it 
 is called box dimension).  Quoting from their paper, a dimension 
 ``corresponds to the possibility of characterizing the 
 ``massiveness'' of sets in metric spaces by the help of the order of 
 growth of the number of elements of their most economical 
 $\eps$-coverings, as $\eps\to0$''.  Set functions retaining the 
 general properties of a dimension (cf.  Theorem~\ref{Thm:Kdim}) have 
 been studied by several authors.  Our choice of the Kolmogorov 
 dimension is due to the fact that it is suitable for the kind of 
 generalization we need in this paper, namely it quite naturally 
 produces a definition of asymptotic dimension.

 In the following, unless otherwise specified, $(X,\d)$ will denote a 
 metric space, $B_X(x,R)$ the open ball in $X$ with centre $x$ and 
 radius $R$, $n_r(\O)$ the least number of open balls of radius $r$ 
 which cover $\O\subset X$, and $\n_r(\O)$ the largest number of 
 disjoint open balls of radius $r$ centered in $\O$.

 The following lemma is proved in \cite{KT}.  Due to some notational 
 difference, we include a proof.

 \begin{Lemma}\label{1.1.1} 
	 $n_r(\O) \geq \n_r(\O)\geq n_{2r}(\O)$. 
 \end{Lemma}
 \begin{proof} 
	 We have only to prove the second inequality when $\n_r$ is 
	 finite.  Let us assume that $\{B(x_i,r)\}_{i=1}^{\n_r(\O)}$ are 
	 disjoint balls centered in $\O$ and observe that, for any $y\in 
	 \O$, $\d(y,\bigcup_{i=1}^{\n_r(\O)} B(x_i,r))<r$, otherwise 
	 $B(y,r)$ would be disjoint from $\bigcup_{i=1}^{\n_r(\O)} 
	 B(x_i,r)$, contradicting the maximality of $\n_r$.  So for all 
	 $y\in \O$ there is $j$ s.t. $\d(y, B(x_j,r))<r$, that is 
	 $\O\subset \bigcup_{i=1}^{\n_r(\O)} B(x_i,2r)$, which implies the 
	 thesis.
 \end{proof}

 Kolmogorov and Tihomirov \cite{KT} defined a dimension for totally 
 bounded metric spaces $X$ as
 \begin{equation}\label{TKdim}
	 \md{X}:=\limsup_{r\to0} \frac{\log n_r(X)}{\log(1/r)}.
 \end{equation}
 Then we may give the following definition.

 \begin{Dfn}\label{Dfn:Kdim}  
	 Let $(X,\d)$ be a metric space.  Then, denoting by $\cb(X)$ the 
	 family of bounded subsets of $X$, the metric Kolmogorov-Tihomirov 
	 dimension of $X$ is
	 $$
	 \md{X}:=\sup_{B\in\cb(X)}\limsup_{r\to0} \frac{\log 
	 n_r(B)}{\log(1/r)}.
	 $$
 \end{Dfn}

 Then the following proposition trivially holds.

 \begin{Prop}\label{1.1.3} 
	 If $\{B_n\}$ is an exhaustion of $X$ by bounded subsets, namely 
	 $B_n$ is increasing and for any bounded $B$ there exists $n$ such 
	 that $B\subseteq B_n$, one has $\md{X}=\lim_n\md{B_n}$.  In 
	 particular,
	 \begin{equation}\label{Kdimformula}
		 \md{X}=\lim_{R\to\infty}\limsup_{r\to0} \frac{\log 
		 n_r(B_X(x,R))}{\log(1/r)}
	 \end{equation}
 \end{Prop}

 \begin{rem} 
	 If bounded subsets of $X$ are not totally bounded, we could 
	 define $\md{X}$ as the supremum over totally bounded subsets.  
	 These two definitions, which agree e.g. on proper spaces, may be 
	 different in general.  For example an orthonormal basis in an 
	 infinite dimensional Hilbert space has infinite dimension 
	 according to Definition~\ref{Dfn:Kdim}, but has zero dimension in 
	 the other case.  A definition of metric dimension which coincides 
	 with $d_0$ on bounded subsets of $\br^p$ has been given by Tricot 
	 \cite{Tricot} in terms of rarefaction indices.
 \end{rem}

 Let us now show that this set function satisfies the basic properties 
 of a dimension \cite{Pontriagin,Tricot}.

 \begin{Thm}\label{Thm:Kdim} 
	 The set function $d_0$ is a dimension, namely it satisfies 
	 \item{$(i)$} If $X\subset Y$ then $\md{X}\leq \md{Y}$.  
	 \item{$(ii)$} If $X_1,X_2\subset X$ then $\md{X_1\cup X_2} = \max 
	 \{ \md{X_1}, \md{X_2} \}$.  
	 \item{$(iii)$} If $X$ and $Y$ are metric spaces, then 
	 $\md{X\times Y} \leq \md{X}+\md{Y}$.
 \end{Thm}
 \begin{proof}
	 Property $(i)$ easily follows from formula~(\ref{Kdimformula}).  
	 \\
	 Now we prove $(ii)$.  The inequality $\md{X_1\cup X_2}$ $\geq$ 
	 $\max \{ \md{X_1},\md{X_2} \}$ follows from $(i)$.  For the 
	 converse inequality, let $x_{i}\in X_{i}$, and set 
	 $\d:=\d(x_{1},x_{2})$, $d_1=\md{X_1}$, $d_2=\md{X_2}$, with e.g. 
	 $d_1\geq d_2$.  If $d_1=\infty$ the property is trivial, so we 
	 may suppose $d_1\in\br$.  Then
	 $$
	 B_{X_1\cup X_2}(x_1,R)\subset B_{X_1}(x_1,R)\cup B_{X_2}(x_2,R+\d)
	 $$
	 therefore
	 \begin{equation}\label{ineq1}
		 n_r( B_{X_1\cup X_2}(x_1,R) ) \leq n_r( B_{X_1}(x_1,R) ) + 
		 n_r( B_{X_2}(x_2,R+\d) ).
	 \end{equation}
	 Now, $\forall R>0$, 
	 $$
	 \limsup_{r\to0} \frac {\log n_r(B_{X_i}(x_{i},R))}{\log 
	 (1/r)}\leq d_i
	 $$ 
	 i.e. $\forall R,\eps>0$ there is $r_0=r_0(\eps,R)$ s.t. for all 
	 $0<r<r_0$, $n_r(B_{X_1}(x_{1},R)) \leq r^{-(d_1+\eps)}$, and 
	 $n_r(B_{X_2}(x_{2},R+\d)) \leq r^{-(d_2+\eps)}$ hence, by 
	 (\ref{ineq1}),
	 $$
	 n_r(B_{X_1\cup X_2}(x,R)) \leq r^{-(d_1+\eps)}(1+r^{d_1-d_2}).
	 $$
	 Finally,
	 $$
	 \limR\limsup_{r\to0} \frac {\log n_r(B_{X_1\cup X_2}(x,R))} {\log 
	 (1/r)} \leq d_1+\eps,
	 $$
	 that is
	 $$
	 \md{X_1\cup X_2}\leq \max\{\md{X_1},\md{X_2}\}+\eps
	 $$
	 and the thesis follows by the arbitrariness of $\eps$.  \\
	 The proof of part $(iii)$ is postponed.
 \end{proof}
 
 Kolmogorov dimension is invariant under bi-Lipschitz maps (also 
 called quasi isometries by some authors), as next proposition shows.
 
 \begin{Prop}\label{1.1.8} 
	 Let $X,Y$ be metric spaces, and $f:X\to Y$ a surjective 
	 bi-Lipschitz map, namely $f$ satisfies
	 $$
	 c_1 \d_X(x_1,x_2) \leq \d_Y(f(x_1),f(x_2)) \leq c_2 \d_X(x_1,x_2).
	 $$
	 Then $\md{X}=\md{Y}$.
 \end{Prop}
 \begin{proof}
	 By hypothesis we have $f(B_X(x,\r/c_2))\subset B_Y(f(x),\r) 
	 \subset f(B_X(x,\r/c_1))$.  So that, with $y_j=f(x_j)$, 
	 $n:=n_r(B_Y(f(x),R))$,
	 \begin{align*}
		 f(B_X(x,R/c_2)) &\subset B_Y(f(x),R) \subset 
		 \bigcup_{j=1}^{n} B_Y(y_j,r) \\
		 &\subset \bigcup_{j=1}^{n} f(B_X(x_j,r/c_1)) = 
		 f(\bigcup_{j=1}^{n} B_X(x_j,r/c_1))
	 \end{align*}
	 which implies $n_{r/c_1}(B_X(x,R/c_2))\leq n_r(B_Y(f(x),R))$.  \\
	 Since bi-Lipschitz maps are injective, we may repeat the same 
	 argument for $f^{-1}$, and we get $n_{c_2r}(B_Y(f(x),c_1R)) \leq 
	 n_r(B_X(x,R))$, so that $n_{r/c_1}(B_X(x,R/c_2))$ $\leq$ 
	 $n_r(B_Y(f(x),R))$ $\leq$ $n_{r/c_2}(B_X(x,R/c_1))$.  Finally
	 \begin{align*}
		 \limsup_{r\to0}\frac {\log n_{r/c_1}(B_X(x,R/c_2)) }{ 
		 \log(c_1/r) - \log c_1} & \leq \limsup_{r\to0}\frac {\log 
		 n_r(B_Y(f(x),R)) }{ \log(1/r)} \\
		 & \leq \limsup_{r\to0}\frac {\log n_{r/c_2}(B_X(x,R/c_1)) }{ 
		 \log(c_2/r) - \log c_2}
	 \end{align*} 
	 and the thesis follows.
 \end{proof}

 \noindent {\it Proof of Theorem~\ref{Thm:Kdim}} (continued). 
	 By the preceding Proposition, we may endow $X\times Y$ with any 
	 metric bi-Lipschitz related to the product metric, i.e.
	 $$
	 \d_{X\times Y}((x_1,y_1), (x_2,y_2))=\max\{ 
	 \d_X(x_1,x_2),\d_Y(y_1,y_2)\}.
	 $$ 
	 Then, by $n_r(B_{X\times Y}((x,y),R))\leq n_r(B_X(x,R))\ 
	 n_r(B_Y(y,R))$, the thesis follows easily.  
 \qed\medskip

 \begin{rem}\label{1.1.9} 
	 Kolmogorov and Tihomirov assign a metric dimension to a totally 
	 bounded metric space $X$ when $\exists\lim_{r\to\infty}$ in 
	 equation (\ref{TKdim}), and consider upper and lower metric 
	 dimensions in the general case.  We observe that if the $\liminf$ 
	 is considered, the classical dimensional inequality 
	 \cite{Pontriagin} stated in Theorem~\ref{Thm:Kdim} (iii) is 
	 replaced by $\md{X\times Y} \geq \md{X}+\md{Y}$.
 \end{rem}

 \subsection{Asymptotic dimension}\label{subsec:asympdim}

 The function introduced in the previous subsection can be used to 
 study local properties of metric spaces.  In this paper we are mainly 
 interested in the investigation of the large scale behavior of these 
 spaces, so we need a different tool.  Looking at 
 equation~(\ref{Kdimformula}), it is natural to set the following

 \begin{Dfn}\label{1.2.1}  
	 Let $(X,\d)$ be a metric space.  We call
	 $$
	 \ad{X}:=\limr\lsup \frac{\log n_r(B_X(x,R)) }{ \log R},
	 $$
	 the {\it asymptotic dimension} of $X$. 
 \end{Dfn}

 Let us remark that, as $n_r(B_X(x,R))$ is nonincreasing in $r$, the 
 function
 $$
 r\mapsto\lsup\frac{\log n_r(B_X(x,R)) }{ \log R}
 $$ 
 is nonincreasing too, so the $\lim_{r\to0}$ exists.

 \begin{Prop}\label{1.2.3} 
	 $\ad{X}$ does not depend on $x$.
 \end{Prop} 
 \begin{proof}
	 Let $x,y\in X$, and set $\d:=\d(x,y)$, so that $B(x,R)\subset 
	 B(y,R+\d) \subset B(x,R+2\d)$.  This implies,
	 \begin{align*}
		 \frac{\log n_r(B(x,R)) }{ \log R} &\leq \frac{\log 
		 n_r(B(y,R+\d)) }{ \log (R+\d)}\ \frac{\log (R+\d) }{ \log R} 
		 \\
		 &\leq \frac{\log n_r(B(x,R+2\d)) }{ \log (R+2\d)}\ \frac{\log 
		 (R+2\d) }{ \log R}
	 \end{align*}
	 so that, taking $\lsup$ and then $\limr$ we get the thesis.  
 \end{proof}

 \begin{Lemma}\label{1.2.4}
	 $$
	 \ad{X}=\limr\lsup \frac{\log \n_r(B_X(x,R)) }{ \log R}
	 $$ 
 \end{Lemma}
 \begin{proof} 
	 Follows easily from lemma \ref{1.1.1}.
 \end{proof}

 \begin{Thm}\label{Thm:adim} 
	 The set function $d_\infty$ is a dimension, namely it satisfies 
	 \item{$(i)$} If $X\subset Y$ then $\ad{X}\leq \ad{Y}$.  
	 \item{$(ii)$} If $X_1,X_2\subset X$ then $\ad{X_1\cup X_2} = \max 
	 \{ \ad{X_1}, \ad{X_2} \}$.  
	 \item{$(iii)$} If $X$ and $Y$ are metric spaces, then 
	 $\ad{X\times Y} \leq \ad{X}+\ad{Y}$.
 \end{Thm}
 \begin{proof}
	 $(i)$ Let $x\in X$, then $B_X(x,R)\subset B_Y(x,R)$ and the claim 
	 follows easily.  \\
	 $(ii)$ By part $(i)$, we get $\ad{X_1\cup X_2} \geq \max \{ 
	 \ad{X_1},\ad{X_2} \}$.  Let us prove the converse inequality.  \\
	 Let $x_i\in X_i$, $i=1,2$, and set $\d=\d(x_1,x_2)$, 
	 $a=\ad{X_1}$, $b=\ad{X_2}$, with e.g. $a\leq b$.  Then, $\forall 
	 \eps,r>0$ $\exists R_0=R_0(\eps,r)$ s.t. $\forall R>R_0$
	 \begin{align*}
		 n_r(B_{X_1}(x_1,R))  & \leq R^{a+\eps} \\
		 n_r(B_{X_2}(x_2,R+\d)) & \leq R^{b+\eps},
	 \end{align*}
	 hence, by inequality (\ref{ineq1}),
	 \begin{align*}
		 n_r(B_{X_1\cup X_2}(x_1,R)) &\leq R^{a+\eps}+R^{b+\eps}\\
		 &= R^{b+\eps}(1+R^{a-b}).\\
	 \end{align*}
	 Finally,
	 $$
	 \frac {\log n_r(B_{X_1\cup X_2}(x_1,R))} {\log R} \leq b+\eps+ 
	 \frac {\log(1+R^{a-b})}{\log R}.
	 $$
	 Taking the $\lsup$ and then the $\limr$ we get
	 $$
	 \ad{X_1\cup X_2}\leq \max\{\ad{X_1},\ad{X_2}\}+\eps
	 $$
	 and the thesis follows by the arbitrariness of $\eps$.  \\
	 The proof of part $(iii)$ is analogous to that of part $(iii)$ in 
	 Theorem~\ref{Thm:Kdim}, where we may use Proposition~\ref{1.2.10} 
	 because bi-Lipschitz maps are rough isometries.  
 \end{proof}

 \begin{rem}\label{1.2.7} 
	 \itm{a} In part $(ii)$ of the previous theorem we considered $X_1$ and 
	 $X_2$ as metric subspaces of $X$.  If $X$ is a Riemannian 
	 manifold and we endow the submanifolds $X_1$, $X_2$ with their 
	 geodesic metrics this property does not hold in general.  A 
	 simple example is the following.  Let $f(t) := (t\cos t,t\sin 
	 t)$, $g(t) := (-t\cos t, -t\sin t)$, $t\geq0$ planar curves, and 
	 set $X,\ Y$ for the closure in $\br^2$ of the two connected 
	 components of $\br^2\setminus (G_f\cup G_g)$, where $G_f,\ G_g$ 
	 are the graphs of $f,\ g$, and endow $X,\ Y$ with the geodesic 
	 metric.  Then $X$ and $Y$ are roughly-isometric to $[0,\infty)$ 
	 (see below) so that $\ad{X}=\ad{Y}=1$, while $\ad{X\cup Y}=2$.\\
	 \itm{b} As for the local case, the choice of the $\limsup$ in 
	 Definition~\ref{1.2.1} is the only one compatible with the 
	 classical dimensional inequality stated in Theorem~\ref{Thm:adim} 
	 $(iii)$.  This, together with the singular traceability property 
	 \ref{thm:ap}, motivates our choice of the $\limsup$ in 
	 Definition \ref{def:NS-inv} for the Novikov-Shubin invariants.
 \end{rem} 

 \begin{Dfn}\label{1.2.8} 
	 Let $X,Y$ be metric spaces, $f:X\to Y$ is said to be a rough 
	 isometry if there are $a\geq1$, $b,\eps\geq0$ s.t. 
	 \itm{i} $a^{-1}\d_X(x_1,x_2)-b \leq \d_Y(f(x_1),f(x_2)) \leq a 
	 \d_X(x_1,x_2)+b$, for all $x_1,x_2\in X$, 
	 \itm{ii} $\bigcup_{x\in X} B_Y(f(x),\eps) = Y$ 
 \end{Dfn}

 It is clear that the notion of rough isometry is weaker then the 
 notion of bi-Lipschitz map introduced in the preceding subsection 
 and, since any compact set is roughly isometric to a point, $d_0$ is 
 not rough-isometry invariant.  We shall show that the asymptotic 
 dimension is indeed invariant under rough isometries.

 \begin{Lemma}\label{1.2.9} {\rm (\cite{Chavel}, Proposition 4.3)} 
	 If $f:X\to Y$ is a rough isometry, there is a rough isometry 
	 $f^-:Y\to X$, with constants $a,b^-,\eps^-$, s.t. 
	 \itm{i} $\d_X(f^-\circ f(x),x)<c_X$, $x\in X$, 
	 \itm{ii} $\d_Y(f\circ f^-(y),y)<c_Y$, $y\in Y$. 
 \end{Lemma}

 \begin{Prop}\label{1.2.10} 
	 Let $X,Y$ be metric spaces, and $f:X\to Y$ a rough isometry.  
	 Then $\ad{X}=\ad{Y}$.
 \end{Prop}
 \begin{proof}
	 Let $x_0\in X$, then for all $x\in B_X(x_0,r)$ we have
	 $$
	 \d_Y(f(x),f(x_0))\leq a \d_X(x,x_0)+b<ar+b
	 $$ 
	 so that
	 $$
	 f(B_X(x_0,r))\subset B_Y(f(x_0),ar+b).
	 $$ 
	 Then, with $n:=n_r(B_Y(f(x_0),aR+b))$, 
	 $$
	 f(B_X(x_0,R)) \subset \bigcup_{j=1}^{n} B_Y(y_j,r),
	 $$
	 which implies 
	 \begin{align*}
		 f^-\circ f(B_X(x_0,R)) & \subset \bigcup_{j=1}^{n} 
		 f^-(B_Y(y_j,r)) \\
		 & \subset \bigcup_{j=1}^{n} B_X(f^-(y_j),ar+b^-).
	 \end{align*}
	 Let $x\in B_X(x_0,R)$, and $j$ be s.t. $f^-\circ f(x)\in 
	 B_X(f^-(y_j),ar+b^-)$, then
	 $$
	 \d_X(x,f^-(y_j))\leq \d_X(x,f^-\circ f(x))+\d_X(f^-\circ 
	 f(x),f^-(y_j))<c_X+ar+b^-,
	 $$ 
	 so that  
	 $$
	 B_X(x_0,R)\subset \bigcup_{j=1}^{n} B_X(f^-(y_j),ar+b^-+c_X),
	 $$
	 which implies $n_{ar+b^-+c_X}(B_X(x_0,R))\leq 
	 n_r(B_Y(f(x_0),aR+b))$.  \\
	 Finally 
	 \begin{align*}
		 \ad{X} & = \limr\lsup\frac {\log n_r(B_X(x_0,R))}{\log R} \\
		 & = \limr\lsup\frac {\log n_{ar+b^-+c_X}(B_X(x_0,R))}{\log R}\\
		 & \leq \limr\lsup\frac {\log n_r(B_Y(f(x_0),aR+b))}{\log R} \\
		 & = \limr\lsup\frac {\log n_r(B_Y(f(x_0),R))}{\log R} \\
		 & = \ad{Y}
	 \end{align*}
	 and exchanging the roles of $X$ and $Y$ we get the thesis.
 \end{proof}
 
 In what follows we show that when $X$ is equipped with a suitable 
 measure, the asymptotic dimension may be recovered in terms of the 
 volume asymptotics for balls of increasing radius, like the local 
 dimension detects the volume asymptotics for balls of infinitesimal 
 radius.

 \begin{Dfn}\label{1.2.12} 
	 A Borel measure $\m$ on $(X,\d)$ is said to be uniformly bounded 
	 if there are functions $\b_1,\b_2$, s.t. $0<\b_1(r)\leq 
	 \m(B(x,r)) \leq \b_2(r)$, for all $x\in X$, $r>0$.  \\
	 That is $\b_1(r):= \inf_{x\in X} \m(B(x,r)) >0$, and $\b_2(r) := 
	 \sup_{x\in X} \m(B(x,r)) <\infty$.
 \end{Dfn}

 \begin{Prop}\label{1.2.13} 
	 If $(X,\d)$ has a uniformly bounded measure, then every ball in 
	 $X$ is totally bounded (so that if $X$ is complete it is locally 
	 compact).
 \end{Prop}
 \begin{proof}
	 Indeed, if there is a ball $B=B(x,R)$ which is not totally 
	 bounded, then there is $r>0$ s.t. every $r$-net in $B$ is 
	 infinite, so $n_r(B)$ is infinite, and $\n_r(B)$ is infinite too.  
	 So that $\b_2(R)\geq \m(B) \geq \sum_{i=1}^{\n_r(B)} \m(B(x_i,r)) 
	 \geq \b_1(r)\n_r(B) = \infty$, which is absurd.  
 \end{proof}

 \begin{Prop}\label{1.2.14} 
	 If $\m$ is a uniformly bounded Borel measure on $X$ then
	 $$
	 \ad{X}=\lsup\frac{\log\m(B(x,R))}{ \log R} .
	 $$ 
 \end{Prop}
 \begin{proof}
	 As $\bigcup_{i=1}^{\n_r(B(x,R))} B(x_i,r) \subset B(x,R+r) 
	 \subset \bigcup_{j=1}^{n_r(B(x,R+r))} B(y_j,r)$, we get
	 $$
	 \b_2(r)n_r(B(x,R+r))\geq \m(B(x,R+r)) \geq \b_1(r)\n_r(B(x,R)) 
	 \geq \b_1(r)n_{2r}(B(x,R)),
	 $$
	 by Lemma \ref{1.1.1}. So that 
	 $$
	 \b_1(r/2)\leq \frac{\m(B(x,R+r/2))}{ n_r(B(x,R))}, \qquad 
	 \frac{\m(B(x,R))}{ n_r(B(x,R))} \leq \b_2(r),
	 $$
	 and the thesis follows easily.
 \end{proof}

 Let us now consider the particular case of complete Riemannian 
 manifolds.
 
 \begin{Prop}\label{2.1.3}
	 Let $M,N$ be complete Riemannian manifolds.  
	 \itm{i} If $M$ is non-compact, then $\ad{M}\geq 1$
	 \itm{ii} If $M$ has bounded geometry, then $\ad{M} = \limR 
	 \frac{\log V(x,R)}{\log R}$
	 \itm{iii} If $M,N$ have bounded geometry, and
	 admit asymptotic dimension in a strong sense, that is
	 $$
	 \ad{M} = \lim_{R\to\infty} \frac{\log V(x,R)}{\log R},
	 $$
	 and analogously for $N$, then 
	 $$
	 \ad{M\times N}=\ad{M}+\ad{N}.
	 $$
 \end{Prop}
 \begin{proof}
	 $(i)$ It follows from Theorem \ref{Thm:adim} $(i)$, and the fact 
	 that there is inside $M$ an unbounded geodesic. \\
	 $(ii)$ It follows from Lemma \ref{2.1.2} that the volume is a 
	 uniformly bounded measure. Therefore the result follows from 
	 Proposition~\ref{1.2.14}. \\
	 $(iii)$ As $vol(B_{M\times N}((x,y),R))=vol(B_M(x,R)) 
	 vol(B_N(y,R))$, we get
	 \begin{align*}
		 \ad{M\times N} &= \limR \frac{\log vol (B_{M\times 
		 N}((x,y),R))}{ \log R} \\
		 & = \limR \frac{\log vol (B_M(x,R))}{ \log R} + \limR 
		 \frac{\log vol (B_N(y,R))}{ \log R} \\
		 & = \ad{M}+\ad{N}.
	 \end{align*}
 \end{proof}

 \begin{rem}
	 Conditions under which the inequality in Theorem~\ref{Thm:adim} 
	 $(iii)$ becomes an equality are often studied in the case of 
	 (local) dimension theory (cf.  \cite{Pontriagin,Salli}).  The 
	 previous Proposition gives such a condition for the asymptotic 
	 dimension.
 \end{rem}

 As the asymptotic dimension is invariant under rough isometries, it 
 is natural to substitute the continuous space with a coarse graining, 
 which destroys the local structure, but preserves the large scale 
 structure.  To state it more precisely, recall (\cite{Chavel}, p.  
 194) that a discretization of a metric space $M$ is a graph $G$ 
 determined by an $\eps$-separated subset $\cg$ of $M$ for which there 
 is a $R>0$ s.t. $M=\cup_{x\in\cg} B_M(x,R)$.  The graph structure on 
 $\cg$ is determined by one oriented edge from any $x\in\cg$ to any 
 $y\in\cg$, $y\neq x$, denoted $<x,y>$, precisely when $\d_M(x,y)<2R$.  
 Define the combinatorial metric on $G$ by $\d_c(x,y):=\inf \{ 
 \sum_{i=0}^n \d(x_i,x_{i+1}): (x_0,\ldots,x_{n+1})\in Path_n(x,y),\ 
 n\in\bn\}$, where $Path_n(x,y) :=\{ (x_0,\ldots,x_{n+1}) : 
 x_i\in\cg,\ x_0=x,\ x_{n+1}=y, <x_i,x_{i+1}>\in G\}$.

 \begin{Prop}\label{2.1.4} {\rm (\cite{Chavel}, Theorem 4.9)} 
	 Let $M$ be a complete Riemannian manifold with Ricci curvature 
	 bounded from below.  Then $M$ is roughly isometric to any of its 
	 discretizations, endowed with the combinatorial metric.  
	 Therefore $M$ has the same asymptotic dimension of any of its 
	 discretizations.
 \end{Prop}

 The previous result, together with the invariance of the asymptotic 
 dimension under rough isometries, shows that, when $M$ has a discrete 
 group of isometries $\Gamma$ with a compact quotient, the asymptotic 
 dimension of the manifold coincides with the asymptotic dimension of 
 the group, hence with its growth (cf.  \cite{GI5}), hence, by the 
 result of Varopoulos \cite{Varopoulos1}, it coincides with the 0-th 
 Novikov-Shubin invariant.  We will generalise this result in 
 subsection \ref{subsec:NSinvariant}.\\

 Let us conclude this subsection with some examples.

 \begin{exmp}\label{1.2.16} 
	 \itm{i} $\br^n$ has asymptotic dimension $n$. 
	 \itm{ii} Set $X:= \cup_{n\in\bz}\{(x,y)\in\br^2 : 
	 \d((x,y),(n,0))<\frac{1}{4} \}$, endowed with the Euclidean 
	 metric, then $\md{X}=2$, $\ad{X}=1$.  
	 \itm{iii} Set $X=\bz$ with the counting measure, then $\md{X}=0$, 
	 and $\ad{X}=1$.  
	 \itm{iv} Let $X$ be the unit ball in an infinite dimensional 
	 Banach space.  Then $\md{X}=+\infty$ while $\ad{X}=0$.
 \end{exmp}

 \begin{exmp}\label{1.2.17} 
	 Set $X:=\{(x,y)\in\br^2: x\geq0, |y|\leq x^\a\}$, endowed with 
	 the Euclidean metric, where $\a\in(0,1]$.  Then $\ad{X}=\a+1$.
 \end{exmp}
 \begin{proof}
	 This metric space has a uniformly bounded Borel measure, the 
	 Lebesgue area, so we can use Proposition \ref{1.2.14}.  Set 
	 $x_0:=(0,0)$, and $B_R:= B_X(x_0,R)$.  Then, if $R\geq 
	 \sqrt{4^{1+1/\a}r^{2/\a}+4^{1+\a}r^2}$, $B_{R}\subset Q_{1}\cup 
	 Q_{2}$, where $Q_1:= \{ (x,y)\in\br^2 : -2r\leq x \leq 
	 (2r)^{1/\a}+2r,\ |y|\leq 4r \}$, and $Q_2:=\{ (x,y)\in\br^2 : 
	 (2r)^{1/\a} \leq x \leq R, \ |y|\leq 2x^\a \}$.  Now, if $x_R>0$ 
	 is s.t. $x_R^2+x_R^{2\a}=R^2$, we get
	 \begin{align*}
		 area(B_R) & \geq \frac2{\a+1} x_R^{\a+1} \\
		 area(Q_1) & = 4r(4r+(2r)^{1/\a}) \\
		 area(Q_2) & = \frac4{\a+1} (R^{1+\a}-(2r)^{1+1/\a}),
	 \end{align*}
	 so that
	 \begin{align*}
		 \limR \frac {(\a+1)\log x_R}{\log R} & \leq \linf \frac {\log 
		 area(B_R)}{\log R} \\
		 & \leq \lsup \frac {\log area(B_R)}{\log R} \leq \a+1
	 \end{align*}
	 and, as $\limR \frac {\log x_R}{\log R} = \lim_{x\to\infty} \frac 
	 {\log x}{\log\sqrt{x^2+x^{2\a}}}=1$, we get the thesis.  
 \end{proof}

 \subsection{Asymptotic dimension of some cylindrical ends}
 \label{subsec:cylindrical}

 In this subsection we want to compare our work with a recent work of 
 Davies'.  In \cite{Davies} he defines the asymptotic dimension 
 of cylindrical ends of a Riemannian manifold $M$ as follows.  Let 
 $E\subset M$ be homeomorphic to $(1,\infty)\times A$, where $A$ is a 
 compact Riemannian manifold.  Set $\dE:=\{1\}\times A$, $E_r:=\{ x\in 
 E: \d(x,\dE) < r\}$, where $\d$ is the restriction of the metric in 
 $M$.  Then $E$ has asymptotic dimension $D$ if there is a positive 
 constant $c$ s.t.
 \begin{equation}\label{e:davies}
	 c^{-1} r^D \leq vol(E_r) \leq cr^D,
 \end{equation}
 for all $r\geq1$.  He does not assume bounded geometry for $E$.  If 
 one does, the two definitions coincide as in the following

 \begin{Prop}\label{2.3.1} 
	 With the above notation, if the volume form on $E$ is a uniformly 
	 bounded measure (as in Definition \ref{1.2.12}), or in particular 
	 if $E$ has bounded geometry (as in Definition \ref{2.1.1}), and 
	 there is $D$ as in {\rm (\ref{e:davies})}, then $\ad{E}=D$.
 \end{Prop}
 \begin{proof}
	 Choose $o\in E$, and set $\d:=\d(o,\dE)$, $\D:=diam(\dE)$.  Then 
	 it is easy to prove that $E_{R-\d-\D} \subset B_E(o,R) \subset 
	 E_{R+\d}$.  \\
	 Then $c^{-1}(R-\d-\D)^D \leq vol(B_E(o,R)) \leq c(R+\d)^D$, and 
	 from \ref{1.2.14} the thesis follows.  
 \end{proof}

 Motivated by (\cite{Davies}, example 16), let us set the following

 \begin{Dfn}\label{2.3.2} 
	 $E$ is a standard end of local dimension $N$ if it is 
	 homeomorphic to $(1,\infty)\times A$, endowed with the metric 
	 $ds^2=dx^2+f(x)^2d\om^2$, and with the volume form 
	 $dvol=f(x)^{N-1}dxd\om$, where $(A,\om)$ is an $(N-1)$-dimensional 
	 compact Riemannian manifold, and $f$ is an increasing smooth 
	 function.
 \end{Dfn}

 \begin{Prop}\label{2.3.3} 
	 The volume form on a standard end $E$ is a uniformly bound\-ed 
	 measure.  Therefore, if $E$ satisfies {\rm (\ref{e:davies})}, we 
	 get $\ad{E}=D$.
 \end{Prop}
 \begin{proof}
	 It is easy to show that, for $(x_0,p_0)\in E$, $r<x_0-1$, \\
	 \begin{align*}
		 [x_0-r/2,x_0+r/2]&\times B_A\left(p_0, 
		 \frac{r/2}{f(x_0+r/2)}\right) \subset B_E((x_0,p_0),r) \\
		 & \subset [x_0-r,x_0+r]\times B_A\left(p_0, 
		 \frac{r}{f(x_0-r)}\right)
	 \end{align*}
	 So that 
	 \begin{align*}
		 \int_{x_0-r/2}^{x_0+r/2}f(x)^{N-1}dx\ & V_A\left(p_0, 
		 \frac{r/2}{f(x_0+r/2)}\right) \leq V_E((x_0,p_0),r) \\
		 &\leq \int_{x_0-r}^{x_0+r}f(x)^{N-1}dx\ V_A\left(p_0, 
		 \frac{r}{f(x_0-r)}\right)
	 \end{align*}
	 which implies
	 \begin{align*}
		 rf(x_0-r/2)^{N-1}\ V_A\left(p_0,\frac{r/2}{f(x_0+r/2)}\right) 
		 &\leq V_E((x_0,p_0),r) \\
		 &\leq 2rf(x_0+r)^{N-1}\ V_A\left(p_0,\frac{r}{f(x_0-r)}\right)
	 \end{align*}
	 As for $x_0\to\infty$, $V_A(p_0,\frac{r}{ f(x_0-r)})\sim c 
	 \left(\frac{r}{f(x_0-r)}\right)^{N-1}$, and the same holds for 
	 $V_A(p_0,\frac{r/2}{f(x_0+r/2)})$, we get the thesis.  
 \end{proof}

 \begin{Cor}\label{2.3.4}  
	 Let $E$ be the standard end of local dimension $N$ and asymptotic 
	 dimension $D$ in {\rm (\cite{Davies}}, example 16), which is 
	 homeomorphic to $(1,\infty)\times S^{N-1}$, endowed with the 
	 metric $ds^2=dr^2+r^{2(D-1)/(N-1)}d\om^2$, and with the volume 
	 form $dvol=r^{D-1}drd^{N-1}\om$.  Then $\ad{E}=D$.
 \end{Cor}

 \begin{rem}
	Observe that $\ad{M}$ makes sense for any metric space, hence for 
	any cylindrical end, while Davies' asymptotic dimension does not.  
	Indeed let $E := (1,\infty)\times S^1$, endowed with the metric 
	$ds^2=dr^2+f(r)^2d\om^2$, and with the volume form 
	$dvol=f(r)drd\om$, where $f(r):=\frac{d}{ dr}(r^2\log r)$.  Then 
	$\ad{E}=2$, but $vol(E_r)$ does not satisfy one of the 
	inequalities in (\ref{e:davies}).
 \end{rem}
 
 Before closing this section we observe that the notion of standard 
 end allows us to construct an example which shows that we could 
 obtain quite different results if we used $\liminf$ instead of 
 $\limsup$ in the definition of the asymptotic dimension.  It makes 
 use of the following function
 $$
 f(x)=
 \begin{cases}
	 \sqrt{x} &  x\in[1,a_1] \\
	 2+ b_{n-1}+ c_{n-1} + (x-a_{2n-1}) & x\in[a_{2n-1},a_{2n}] \\
	 2+ b_{n-1}+ c_{n} + \sqrt{ x-a_{2n}+1 } & x\in[a_{2n},a_{2n+1}]
 \end{cases}
 $$
 where $a_0:=0$, $a_{n}-a_{n-1}:=2^{2^n}$, $b_{n}:= \sum_{k=1}^{n} 
 \sqrt{ 2^{2^{2k+1}}+1 }$, $c_{n}:= \sum_{k=1}^{n} ( 2^{2^{2k}}-1 )$, 
 $n\geq1$.

 \begin{Prop}\label{2.3.5} 
	 Let $M$ be the Riemannian manifold obtained as a $C^\infty$ 
	 regularization of $C\cup_\f E$, where $C:= \{ (x,y,z)\in\br^3 : 
	 (x-1)^2+y^2+z^2=1,\ x\leq1 \}$, with the Euclidean metric, 
	 $E:=[1,\infty)\times S^1$, endowed with the metric 
	 $ds^2=dx^2+f(x)^2d\om^2$, and with the volume form 
	 $dvol=f(x)dxd\om$, where $\f$ is the identification of $\{ 
	 y^2+z^2=1,\ x=1 \}$ with $\{1\}\times S^1$.  Then the volume form 
	 is a uniformly bounded measure, $\ad{M}\geq2$ but 
	 $\underline{d}_\infty(M)\leq3/2$, where 
	 $\underline{d}_\infty(M):=\limr\liminf_{R\to\infty}\frac{\log 
	 n_r(B_M(x,R)) }{ \log R}$.
 \end{Prop}
 \begin{proof}
	 Set $o:=(0,0,0)\in M$, then it is easy to see that, for 
	 $n\to\infty$, $a_{n}\sim 2^{2^n}$, $b_{n}\sim c_{n} \sim 
	 2^{2^{2n}}$, and
	 \begin{align*}
		 area(B_M(o,a_{2n})) & \sim \frac12 a_{2n}^2 \\
		 area(B_M(o,a_{2n-1}) )& \sim \frac53 a_{2n-1}^{3/2} 
	 \end{align*}
	 so that, calculating the limit of $\frac{\log 
	 area(B_M(o,R))}{\log R}$ on the sequence $R=a_{2n}$ we get $2$, 
	 while on the sequence $R=a_{2n-1}$ we get $3/2$.  The thesis 
	 follows easily, using Proposition \ref{1.2.14}.  
 \end{proof}
 
 \subsection{The asymptotic dimension and the 0-th Novikov Shubin 
 invariant} \label{subsec:NSinvariant}

 In this subsection we show that, for open manifolds of \bg and 
 satisfying an isoperimetric inequality due to Grigor'yan 
 \cite{Grigoryan94}, the asymptotic dimension coincides with the 
 $0$-th Novikov-Shubin invariant. Let us start by recalling a recent 
 result of Coulhon-Grigor'yan which is crucial for the following.

 \begin{Thm}\label{2.1.6} 
	 {\rm (\cite{CG96}, Corollary 7.3) (\cite{Grigoryan94}, 
	 Proposition 5.2)} \\
	 Let $M$ be a complete Riemannian manifold, and set $\l_1(U)$ for 
	 the first Dirichlet eigenvalue of $\D$ in $U$.  Then the 
	 following are equivalent 
	 \itm{i} there are $\a,\ \b>0$ s.t. for all $x\in M$, $r>0$, and 
	 all regions $U\subset B(x,r)$, 
	 \begin{equation}\label{e:isop}
		 \l_1(U) \geq \frac{\a }{ r^2} \left( \frac{V(x,r)}{ vol(U)} 
		 \right)^\b
	 \end{equation}
	 \itm{ii} there are $A,\ C,\ C'>0$ s.t. for all $x\in M$, $r>0$,
	 \begin{align}
		 V(x,2r) & \leq A V(x,r) \label{e:voldouble}\\
		 \frac{C}{ V(x,\sqrt{r})} & \leq H_{0}(r,x,x) \leq 
		 \frac{C'}{V(x,\sqrt{r})} .  \label{e:stime}
	 \end{align}
 \end{Thm}

 Condition (\ref{e:isop}) is introduced in \cite{Grigoryan94} and 
 called isoperimetric inequality, whereas inequality~(\ref{e:voldouble}) 
 is introduced in \cite{CG96} and called the volume doubling property.

 \begin{Cor}\label{2.1.7} 
	 Let $M$ be a complete Riemannian manifold of bounded geometry, 
	 and assume one of the equivalent properties of the previous 
	 Theorem.  \\
	 Then $\ad{M}=\limsup_{t\to\infty} \frac{-2\log 
	 H_{0}(t,x_0,x_0)}{\log t}$
 \end{Cor}
 \begin{proof} 
	 Follows from Theorem \ref{2.1.2} and estimates (\ref{e:stime}).
 \end{proof}
 
 \begin{rem}
	The previous result shows that there are some connections between 
	the asymptotic dimension of a manifold and the notion of dimension 
	at infinity for semigroups (in our case the heat kernel semigroup) 
	considered by Varopoulos (see \cite{VSC}).
 \end{rem}

 The volume doubling property is a weak form of polynomial growth 
 condition, but still guarantees the finiteness of the asymptotic 
 dimension (for manifolds of bounded geometry).

 \begin{Prop}
	 Let $M$ be a complete Riemannian manifold of bounded geometry, 
	 and suppose the volume doubling property (\ref{e:voldouble}) 
	 holds.  Then $M$ has finite asymptotic dimension.
 \end{Prop}
 \begin{proof}
	 Let $R>1$, and $n\in\bn$ be s.t. $2^{n-1}<R \leq 2^n$.  Then 
	 $V(x,R)\leq V(x,2^n)\leq A^n V(x,1)$, so that
	 $$
	 1 \leq \frac{V(x,R) }{ V(x,1)} \leq A^n \leq A R^{\log_2 A}.
	 $$
	 Therefore $\ad{M} = \limsup_{R\to\infty} \frac{\log V(x,R) }{ 
	 \log R} \leq \log_2 A$.  
 \end{proof}

 From now on $M$ is an open manifolds of \bg and satisfying the
 isoperimetric inequality (\ref{e:isop}). Then it has finite asymptotic 
 dimension, which we show to coincide with the 0-th Novikov-Shubin 
 invariant. First we need
 
 \begin{Prop}\label{p:exbyballs}	
	 \itm{i} For any $x,\ y\in M$, $r>0$, if $B(x,r)\cap B(y,r) \neq 
	 \emptyset$, then
	 $$
	 \g^{-1}\leq\frac{V(x,r)}{V(y,r)}\leq\g.
	 $$
	 \itm{ii} There is a sequence $n_{k}\in\bn$ s.t. $\{B(x,n_{k})\}$ 
	 is a regular exhaustion of $M$.
 \end{Prop}
 \begin{proof} 
	 \itm{i} The inequality easily follows by a result of Grigor'yan 
	 (\cite{Grigoryan94}, Proposition 5.2), where it is shown that the 
	 isoperimetric inequality above implies the existence of a 
	 constant $\g$ such that
	 $$
	 \g^{-1}\left(\frac{R}{r}\right)^{\a_{1}} \leq 
	 \frac{V(x,R)}{V(y,r)} \leq \g\left(\frac{R}{r}\right)^{\a_{2}}
	 $$
	 for some positive constants $\a_{1},\ \a_{2}$, for any $R\geq r$, 
	 and $B(x,R)\cap B(y,r)\neq\emptyset$. 
	 \itm{ii} The statement follows from the fact that the volume 
	 doubling property implies subexponential (volume) growth, so that 
	 the result is contained in (\cite{Roe1}, Proposition 6.2).
 \end{proof}

 \begin{Thm}\label{a0=dinfinity} 
	 Let $M$ be an open manifold of \bg and satisfying Grigor'yan's 
	 isoperimetric inequality (\ref{e:isop}), endowed with the regular 
	 exhaustion $\ck$ given by Proposition \ref{p:exbyballs} $(ii)$.  
	 Then the asymptotic dimension of $M$ coincides with the 0-th 
	 Novikov-Shubin invariant, namely $d_\infty(M)=\a_0(M,\ck)$. In 
	 particular $\a_{0}$ is independent of the limit procedure $\om$.
 \end{Thm}
 \begin{proof} 
	 First, from Theorem \ref{2.1.6} and the previous Proposition, we get
	 \begin{align*}
		 \frac{C\g^{-1}}{V(o,\sqrt{t})} &\leq \frac {\int_{B(o,r)} 
		 \frac{C}{V(x,\sqrt{t})} dvol(x)} {V(o,r)} \leq 
		 \frac{\int_{B(o,r)}H_{0}(t,x,x)dvol(x)}{V(o,r)}\\
		 &\leq \frac{\int_{B(o,r)} \frac{C'}{V(x,\sqrt{t})} dvol(x)} 
		 {V(o,r)} \leq \frac{C'\g}{V(o,\sqrt{t})}
	 \end{align*}
	 therefore, by definition of the trace $Tr_{\ck}$,
	 $$
	 \frac{C\g^{-1}}{V(o,\sqrt{t})} \leq Tr_{\ck}(\e{-t\D}) \leq 
	 \frac{C'\g}{V(o,\sqrt{t})}
	 $$
	 hence, finally,
	 \begin{align*}
		 d_\infty(M) & = 
		 2\limsup_{t\to\infty}\frac{\log(V(o,t))}{2\log t} = 
		 2\limsup_{t\to\infty}\frac{\log(C'\g
		 V(o,\sqrt{t})^{-1})}{\log\frac{1}{t}}\\
		 \leq \a'_0(M) & \equiv 2\limsup_{t\to\infty}\frac 
		 {\log\t(\e{-t\D})}{\log\frac{1}{t}} 
		 \leq2\limsup_{t\to\infty}\frac {\log(C\g^{-1}  
		 V(o,\sqrt{t})^{-1})}{\log\frac{1}{t}}\\
		 & = 2\limsup_{t\to\infty}\frac {\log(V(o,t))}{2\log t} = 
		 d_\infty(M).
	 \end{align*}
	 The thesis then follows from Proposition \ref{p:a0=a'0}.
 \end{proof}

 \begin{ack}
	 We would like to thank D.~Burghelea, M.~Farber, L.~Friedlander, 
	 P.~Piazza, M.~Shubin for conversations.  We also thank I.~Chavel, 
	 E.B.~Davies, A.~Grigor'yan, P.~Li, L.~Saloff-Coste for having 
	 suggested useful references on heat kernel estimates.
 \end{ack}

%%%%%%REFERENCES%%%%%%%%


\begin{thebibliography}{99}

 \bibitem{AGPS} S. Albeverio, D.~Guido, A.~Ponosov, S.~Scarlatti.  
 {\it Singular traces and compact operators}.  J. Funct.  Anal., {\bf 
 137} (1996), 281--302.

 \bibitem{Atiyah} M. F. Atiyah.  {\it Elliptic operators, discrete 
 groups and von Neumann algebras}.  Soc.  Math.  de France, 
 Ast\'erisque {\bf 32--33} (1976), 43--72.
 
 \bibitem{BFK} D. Burghelea, L. Friedlander, T. Kappeler.  {\it 
 Relative torsion for representations in finite type Hilbert modules}.  
 Preprint dg-ga/9711018. 

 \bibitem{BE} I. Buttig, J. Eichhorn.  {\it The heat kernel for 
 $p$-forms on manifolds of bounded geometry}.  Acta Sci.  Math., {\bf 
 55} (1991), 33--51.
 
 \bibitem{CFM} A. Carey, M. Farber, V. Mathai.  {\it Determinant lines, von 
 Neumann algebras and L$^{2}$-torsion}. Journal f\" ur die reine und 
 angewandte Mathematik, {\bf 484} (1997), 153--181.

 \bibitem{Chavel} I. Chavel.  {\it Riemannian geometry - A modern 
 introduction}.  Cambridge Univ.  Press, Cambridge, 1993.

 \bibitem{ChavelFeldman} I. Chavel, E. A. Feldman.  {\it Modified 
 isoperimetric constants, and large time heat diffusion in Riemannian 
 manifolds}.  Duke J. Math., {\bf 64} (1991), 473--499.

 \bibitem{Christensen} E. Christensen.  {\it Non commutative 
 integration for monotone sequentially closed C$^{*}$-algebras}.  
 Math.  Scand., {\bf 31} (1972), 171--190.

 \bibitem{Combes} F. Combes. {\it Poids sur une C$^{*}$-alg\`ebre}. 
 J. Math. pures et appl., {\bf 47} (1968), 57--100.
 
 \bibitem{CoCMP} A. Connes.  {\it The action functional in non 
 commutative geometry}.  Commun.  Math.  Phys., {\bf 117} (1988), 
 673--683.

 \bibitem{Co} A. Connes.  {\it Non Commutative Geometry}.  Academic 
 Press, 1994.

 \bibitem{CoMo} A. Connes, H. Moscovici.  {\it The local index formula 
 in noncommutative geometry}.  GAFA, {\bf 5} (1995), 174--243.

 \bibitem{Cou} T. Coulhon.  {\it Dimension \`a l' infini d' un 
 semi-groupe analytique}.  Bull.  Sc.  Math., 2\`eme s\'erie, {\bf 
 114} (1990), 485--500.

 \bibitem{CG96} T. Coulhon, A. Grigor'yan.  {\it On-diagonal lower 
 bounds for heat kernels and Markov chains}.  Duke Math. J., {\bf 89} 
 (1997), 133--199.

 \bibitem{Davies88} E. B. Davies.  {\it Gaussian upper bounds for the 
 heat kernels of some second--order operators on Riemannian 
 manifolds}.  J. Funct.  Anal., {\bf 80} (1988), 16--32.

 \bibitem{DaviesJOT} E. B. Davies.  {\it Pointwise bounds on the space 
 and time derivatives of heat kernels}.  J. Oper.  Th., {\bf 21} 
 (1989), 367--378.

 \bibitem{Daviesbook} E. B. Davies.  {\it Heat kernels and spectral 
 theory}.  Cambridge Univ.  Press, Cambridge, 1989.

 \bibitem{Davies} E. B. Davies.  {\it Non-gaussian aspects of heat 
 kernel behaviour}.  J. London Math. Soc., {\bf 55} (1997), 105--125.

 \bibitem{Dixmier} J. Dixmier.  {\it Existence de traces non 
 normales}.  C.R. Acad.  Sci.  Paris, {\bf 262} (1966), 1107--1108.

 \bibitem{DixmierC} J. Dixmier.  {\it C$^*$-algebras}.  North-Holland 
 Publ., Amsterdam, 1977.
 
 \bibitem{Dodziuk} J. Dodziuk.  {\it De Rham-Hodge theory for 
 $L^2$-cohomology of infinite coverings}.  Topology, {\bf 16} (1977), 
 157--165.

 \bibitem{ES} D. V. Efremov, M. A. Shubin.  {\it Spectrum distribution 
 function and variational principle for automorphic operators on 
 hyperbolic space}.  S\'eminaire Equations aux Deriv\'ees Partielles, 
 Ecole Polytechnique, Centre de Math\'ematique, Expos\'e VII 
 (1988-89).

 \bibitem{FK} T. Fack, H. Kosaki.  {\it Generalized s-numbers of 
 $\t$-measurable operators}.  Pacific J. Math., {\bf 123} (1986), 269.

 \bibitem{Falconer} K. Falconer. {\it Fractal geometry}. J. Wiley, 
 New York, 1990.
 
 \bibitem{Farber} M. Farber.  {\it Geometry of growth: approximation 
 theorems for $L^2$ invariants}.  Math.  Ann., {\bf 311} (1998), 
 335--375.
 
 \bibitem{Grigoryan94} A. Grigor'yan.  {\it Heat kernel upper bounds 
 on a complete non-compact manifold}.  Revista Matematica 
 Iberoamericana, {\bf 10} (1994), 395--452.

 \bibitem{GS} M. Gromov, M. Shubin.  {\it Von Neumann spectra near 
 zero}.  Geometric and Functional Analysis, {\bf 1} (1991), 375--404.

 \bibitem{GI1} D. Guido, T. Isola.  {\it Singular traces for 
 semifinite von~Neumann algebras}.  Journal of Functional Analysis, 
 {\bf 134} (1995), 451--485.

 \bibitem{GI4} D. Guido, T. Isola.  {\it Noncommutative Riemann 
 integration and singular traces for C$^{*}$- algebras}.  Preprint.

 \bibitem{GI5} D. Guido, T. Isola.  {\it Singular traces and 
 Novikov-Shubin invariants}.  Preprint.

 \bibitem{HeSt} E. Hewitt, K. Stromberg. {\it Real and abstract 
 analysis}. Springer, New York, 1975.

 \bibitem{KT} I. Kolmogoroff, M.G. Tihomirov.  {\it $\eps$-Entropy and 
 $\eps$-capacity of sets in functional spaces}.  A. M. S. Transl., 
 {\bf 17} (1961), 277--364.

 \bibitem{Lott} J. Lott.  {\it Heat kernels on covering spaces and 
 topological invariants}.  J. Diff.  Geom., {\bf 35} (1992), 471--510.
 
 \bibitem{Lott1} J. Lott.  {\it Invariant Currents on Limit Sets}. 
 Preprint math/9807025.
 
 \bibitem{Lueck} W. L\" uck.  {\it L$^{2}$-invariants of regular 
 coverings of compact manifolds and CW-complexes}.  To appear in 
 "Handbook of geometry", R.J. Davermann and R.B. Sher (eds.), Elsevier 
 (1997)
 
 \bibitem{Lueck1} W. L\" uck.  {\it Dimension theory of arbitrary 
 modules over finite von Neumann algebras and $L^2$-Betti numbers I: 
 Foundations}.  Journal f\"ur reine und angewandte Mathematik, {\bf 
 495} (1998),135 -- 162.

 \bibitem{Moore-Schochet} C. C. Moore, C. Schochet.  {\it Global 
 analysis on foliated spaces}.  Springer, New York, 1988.

 \bibitem{Ne} E. Nelson.  {\it Notes on non-commutative integration}.  
 Journ.  Funct.  An., {\bf 15} (1974), 103.

 \bibitem{NS1} S. P. Novikov, M. A. Shubin.  {\it Morse theory and von 
 Neumann {\rm II}${}_1$ factors}.  Doklady Akad.  Nauk SSSR, {\bf 289} 
 (1986), 289--292.

 \bibitem{NS2} S. P. Novikov, M. A. Shubin.  {\it Morse theory and von 
 Neumann invariants on non-simply connected manifolds}.  Uspekhi Math.  
 Nauk, {\bf 41}, 5 (1986), 222--223 (in Russian).

 \bibitem{Ped} G. K. Pedersen.  {\it C$^{*}$-algebras and their 
 automorphism groups}.  Academic Press, London, 1979.

 \bibitem{Pontriagin} A. V. Arkhangel'skii, L. S. Pontryagin (eds.).  
 {\it General Topology I. Basic concepts and constructions.  Dimension 
 theory}.  Enc.  Math.  Sci.  {\bf 17}, Springer, New York, 1988.

 \bibitem{QV} J. Quaegebeur, J. Verding.  {\it A construction for 
 weights on C$^{*}$-algebras.  Dual weights on C$^{*}$-crossed 
 products}.  Preprint.
 
 \bibitem{Roe1} J. Roe.  {\it An index theorem on open manifolds.  I}.  
 J. Diff.  Geom., {\bf 27} (1988), 87--113.

 \bibitem{Roe2} J. Roe.  {\it An index theorem on open manifolds.  
 II}.  J. Diff.  Geom., {\bf 27} (1988), 115--136.

 \bibitem{RoeBetti} J. Roe.  {\it On the quasi-isometry invariance of 
 L$^{2}$ Betti numbers}.  Duke Math.  J., {\bf 59} (1989), 765--783.

 \bibitem{RoeCoarse} J. Roe.  {\it Coarse cohomology and index theory 
 on complete Riemannian manifolds}.  Memoirs A.M.S., {\bf 497} (1993), 
 1--90.

 \bibitem{Salli} A. Salli.  {\it On the Minkowski dimension of 
 strongly porous fractal sets in $\br^{n}$}.  Proc.  London Math.  
 Soc., {\bf 62} (1991), 353--372.

 \bibitem{Se} I.E. Segal.  {\it A non-commutative extension of 
 abstract integration}.  Ann.  Math.  {\bf57} (1953) 401.

 \bibitem{Stratila} S. Stratila.  {\it Modular theory in operator 
 algebras}, Abacus Press, Tunbridge Wells, England, 1981.
 
 \bibitem{Tricot} C. Tricot.  {\it Two definitions of fractional 
 dimension}.  Math.  Proc.  Camb.  Phil.  Soc., {\bf 91} (1982), 
 57--74.
 
 \bibitem{Varopoulos1} N. T. Varopoulos.  {\it Random walks and 
 Brownian motion on manifolds}.  Symposia Mathematica, {\bf XXIX} 
 (1987), 97--109.
 
 \bibitem{Varopoulos2} N. T. Varopoulos.  {\it Brownian motion and 
 random walks on manifolds}.  Ann.  Inst.  Fourier, Grenoble, {\bf34}, 
 2, (1984), 243--269.

 \bibitem{VSC} N. T. Varopoulos, L. Saloff-Coste, T. Couhlon.  {\it 
 Analysis and geometry on groups}.  Cambridge Univ.  Press, Cambridge, 
 1992.

 \bibitem{Wo} M. Wodzicki.  {\it Noncommutative residue, part I. 
 Fundamentals}.  In ``{\it K-theory, arithmetic and geometry}'' 
 (Moskow, 1984-86), pp.  320--399, Lecture Notes in Math.  {\bf 1289}, 
 Springer, Berlin, 1987.

 \bibitem{Wolf} J. Wolf.  {\it Essential self-adjointness for the 
 Dirac operator and its square}.  Indiana Univ.  Math.  J. {\bf 22} 
 (1972/73), 611--640.
 
\end{thebibliography}
\end{document}